\documentclass{amsproc}

\newcommand{\bz}{{\mathbb B}}
\newcommand{\cz}{{\mathbb C}}

\newcommand{\nz}{{\mathbb N}}
\newcommand{\rz}{{\mathbb R}}
\newcommand{\gz}{{\mathbb Z}}

\newcommand{\calA}{\mathcal{A}}
\newcommand{\calC}{\mathcal{C}}
\newcommand{\calD}{\mathcal{D}}
\newcommand{\calE}{\mathcal{E}}
\newcommand{\calH}{\mathcal{H}}
\newcommand{\calK}{\mathcal{K}}
\newcommand{\calL}{\mathcal{L}}
\newcommand{\calM}{\mathcal{M}}

\newcommand{\calS}{\mathcal{S}}

\newcommand{\amax}{A_{\text{\rm max}}}
\newcommand{\amin}{A_{\text{\rm min}}}
\newcommand{\ci}{\mathcal{C}^\infty}
\newcommand{\cicomp}{\mathcal{C}^\infty_{\text{\rm comp}}}

\newcommand{\dbar}{d\hspace*{-0.08em}\bar{}\hspace*{0.1em}}
\newcommand{\eps}{\varepsilon}

\newcommand{\hsgpb}{\mathcal{H}^{s,\gamma}_p({\mathbb B})}

\newcommand{\im}{\text{\rm Im}\,}
\newcommand{\intb}{\text{\rm int}\,{\mathbb B}}

\newcommand{\pit}{\,{\widehat{\otimes}}_\pi\,}
\newcommand{\re}{\text{\rm Re}\,}
\newcommand{\rp}{{\mathbb R}_+}

\newcommand{\skp}[2]{\langle#1,#2\rangle}
\newcommand{\smsum}{\mathop{\mbox{\Large$\sum$}}}
\newcommand{\spk}[1]{\left<#1\right>}
\newcommand{\st}{\mbox{\boldmath$\;|\;$\unboldmath}}
\newcommand{\supp}{\mathrm{supp}\;}
\newcommand{\trinorm}[1]%
    {|\hspace*{-1pt}|\hspace*{-1pt}|#1|\hspace*{-1pt}|\hspace*{-1pt}|}

\newcommand{\Mr}{M\!R}

\newcommand{\sfrac}[2]{\frac{#1}{#2}} 

\newtheorem{theorem}{Theorem}[section]
\newtheorem{lemma}[theorem]{Lemma}
\newtheorem{corollary}[theorem]{Corollary}
\newtheorem{proposition}[theorem]{Proposition}
\newtheorem{definition}[theorem]{Definition}
\newtheorem{example}[theorem]{Example}
\newtheorem{remark}[theorem]{Remark}

\numberwithin{equation}{section}

\begin{document}

\title[Cone Operators: Maximal Regularity and Parabolic Equations]
      {Differential Operators on Conic Manifolds:\\
       Maximal Regularity and Parabolic Equations}
\author{S.\ Coriasco$^\dagger$}
\address{Universit\'a di Torino, Dipartimento di Matematica,
         Via Carlo Alberto 10, 10123 Torino, Italy}
\email{coriasco@dm.unito.it}
\author{E.\ Schrohe}
\address{Universit\"at Potsdam, Institut f\"ur Mathematik, 
         Postfach 60 15 53, 14415 Potsdam, Germany}
\email{schrohe@math.uni-potsdam.de}
\author{J.\ Seiler}
\address{Universit\"at Potsdam, Institut f\"ur Mathematik, 
         Postfach 60 15 53, 14415 Potsdam, Germany}
\email{seiler@math.uni-potsdam.de}
\date{\today}
\keywords{Semilinear parabolic equations, Manifolds with 
conical singularities.\newline
\mbox{\hspace{0.65cm}$^\dagger$}Supported by the E.U. 
          Research and Training Network ``Geometric Analysis''.}

\begin{center}
{\Large \bf 
Differential Operators on Conic Manifolds:\\[.5ex]
       Maximal Regularity and Parabolic Equations}

\vspace{2.5ex}
{\large 
S.\ Coriasco\footnote{
Supported by the E.U.\ Research and Training
Network ``Geometric Analysis''.\\ 
\hspace*{.66cm}
\hspace*{.66cm}
},
E.\ Schrohe, and J.\ Seiler}

\vspace{1cm}
{\em \footnotesize Dedicated to the memory of Pascal Laubin}

\vspace{1cm}

\begin{minipage}{12.5cm}
\setlength{\baselineskip}{8pt}
\footnotesize {\sc Abstract.} 
    We study an elliptic differential operator $A$ on a manifold with 
    conic points. Assuming $A$ to be defined on the smooth functions 
    supported away from the singularities, we first address the 
    question of possible closed extensions of $A$ to $L_{p}$ Sobolev 
    spaces and then explain how additional ellipticity conditions 
    ensure maximal regularity for the operator $A$. 
    Investigating the Lipschitz continuity of the maps 
    $f(u)=|u|^\alpha$, $\alpha\ge1$, and $f(u)=u^\alpha$, $\alpha\in\nz$, 
    and using a result of Cl\'ement and Li, 
    we finally show unique solvability of a quasilinear  
    equation of the form $(\partial_{t} - a(u)\Delta) u = f(u)$ in suitable
spaces.
\end{minipage}
\end{center}


\markboth{\uppercase{Cone Operators: Maximal Regularity and Parabolic Equations}}
         {\uppercase{S.\ Coriasco, E.\ Schrohe, J.\ Seiler}}

\tableofcontents

\setlength{\parskip}{5pt}

\section{Introduction}
Parabolic equations and associated initial value problems or boundary 
value problems are common models appearing in science and engineering.
A well-known example is the mixed initial-boundary value
problem for the heat equation
\begin{equation}
    \label{heat}
    \begin{cases}
	\partial_{t} u(t,x) - \Delta u(t,x) = g(t,x) & 
	                      \mbox{on ${]0,T[}\times \Omega$,}
	\\
	u(0,x) = u_{0}(x) & \mbox{on }\Omega, 
	\\
	u(t,x)|_{\partial \Omega} = u_{1}(x) & 
	                       \mbox{for $t \in {]0,T[}$},
     \end{cases}
\end{equation}
where $\Omega$ is a domain (or manifold) with smooth 
boundary $\partial \Omega$.

A typical approach to solve \eqref{heat} consists in rewriting 
it as an abstract evolution equation
\begin{equation}
    \label{abspara}
    \begin{cases}
	\dot{u}(t) + A u(t) = g(t) & 
	                      \mbox{on ${]0,T[},$}
	\\
	u(0) = u_{0}
    \end{cases}
\end{equation}
with an unbounded
operator $A$ on a Banach space $E_{0}$, whose domain $E_{1} = \calD(A)$ 
is densely and continuously embedded into $E_{0}$ and incorporates the choice
of the boundary
condition. The investigation of existence, uniqueness, and regularity of the 
solution to the problem \eqref{abspara}, both in the $\calC^{1}$ and 
$L_{q}$ setting, has attracted the attention of
many authors,  see, e.g., Amann \cite{Aman}, Arendt et al. \cite{Aren},
Da Prato and Grisvard \cite{DaprGr}, Lunardi \cite{Luna}, 
and Pr\"uss \cite{Prus}.
 
\eqref{heat} can be viewed as a special case of the following 
semilinear problem
\begin{equation}
    \label{GinzLan}
    \begin{cases}
	\partial_{t} u(t,x) + A u(t,x) = f(t,u(t,x)) + g(t,x) & 
	                      \mbox{on ${]0,T[} \times \Omega,$}
	\\
	u(0,x) = u_{0}(x) & \mbox{on }\Omega,
	\\
	u(t,x)|_{\partial \Omega} = u_{1}(x) & 
	                       \mbox{for $t \in {]0,T[}$}.
    \end{cases}
\end{equation}
For $A=-\Delta$, the problem \eqref{GinzLan} is a so-called 
reaction-diffusion system which models phenomena in physics, chemistry 
and biology (see, e.g., \cite{CaHa}). A problem in superconductivity 
for example is described by the so-called (non-stationary) 
Ginzburg-Landau equation, where $\Omega$ is a cube
in $\rz^{3}$, $A=-\Delta$ is the Laplacian, $g \equiv 0$, and the 
nonlinearity is $f(t,u)=u-u^{3}$.

More generally, the operator $A$ might also depend on $u$.
In \cite{ClLi}, Cl\'ement and Li developed a method for solving the 
 quasilinear abstract problem 
\begin{equation}
    \label{semi}
    \begin{cases}
	\dot{u}(t) + A(u) u(t) = f(t,u(t)) + g(t) & 
	                      \mbox{on ${]0,T[}$}
	\\
	u(0) = u_{0}
    \end{cases}
\end{equation}
in the $L_{q}$ setting, which relies
on the properties of the linear problem \eqref{abspara} associated with
$A(u_0)$. The main requirement is that $A(u_0)$ be of ``maximal $L_{q}$ 
regularity'', which, in short, means that, for 
every choice of $g \in L_{q}([0,T],E_{0})$ and
$u_{0}$ in the real interpolation space 
$(E_{1},E_{0})_{\frac{1}{q},q}$, 
\eqref{abspara} admits a unique solution $u \in 
L_{q}([0,T],E_{1})\cap W_{q}^{1}([0,T],E_{0})$. 
This together with appropriate Lipschitz continuity of $A$ and $f$ gives
existence and uniqueness of the solution to \eqref{semi}.

In this paper we consider the case where $A$ is a differential 
operator on a manifold with conic singularities which we denote by 
$B$. Formally, $B$ is a compact Hausdorff space which is a smooth 
manifold outside a finite number of singular points, while, near each 
of these points, it has the structure of a cone whose cross-section is 
a smooth closed manifold. 

In order to describe the class of operators 
we treat in this paper, we blow up $B$ at the singularities so that 
we obtain a manifold $\bz$ with boundary $X=\partial \bz$.
When we speak of a differential operator $A$ on the conic manifold $B$ or a 
{\em cone differential operator}, we shall mean a differential operator on $\intb$, 
the interior of $\bz$, which has a Fuchs type degeneracy near 
the boundary, i.e., with respect to a splitting of coordinates 
$(t,x)\in {[0,1[}\times X$ near the boundary, it is of the form
  \begin{equation}\label{fuchs0}
    A=t^{-\mu}\smsum_{j=0}^\mu a_{j}(t)(-t\partial_{t})^j,\qquad
    a_{j}\in\ci([0,1[,\text{\rm Diff}^{\,\mu-j}(X))
  \end{equation}
(note that from now on we shall use  $t$ no longer as the time variable;
instead it corresponds to the distance from the boundary in this neighborhood).

While one should keep the intuitive picture of the conic manifold $B$ 
in mind, it is important and a great simplification that all the 
analysis takes place on $\bz$. In fact, the only way the singularity 
of the underlying manifold then enters into the considerations
is through the particular form of the operators we study. The 
Fuchs type degeneracy encodes that the singularities are conic; other 
types of singularities can be modelled by corresponding degeneracies, 
see, e.g., Schulze \cite{Schu2}, Melrose \cite{Melr}, Mazzeo
\cite{Mazz}.

The choice of Fuchs type operators is motivated by two observations. 
First, consider $\rz^{n+1}$ as the cone over $S^n$ with vertex at the 
origin. The blow-up and the use of variables $(t,x)$ described 
above correspond to the choice of polar coordinates. A simple 
computation shows that every differential operator with smooth 
coefficients on $\rz^{n+1}$ then takes the form \eqref{fuchs0}. 
Note, however, that the class of Fuchs type operators is considerably 
larger. It includes operators with discontinuous coefficients at $0$
(only the radial limits have to exist).
The second observation is that the Laplace-Beltrami operator with respect 
to a Riemannian metric with a conic degeneracy also has this form, cf.\ 
Example \ref{laplace1}.

As the main result of this paper we shall show in Section \ref{semi5} 
that Cl\'ement and Li's method yields $L_{q}$ solvability of certain problems of type 
\eqref{semi}. 
Our argument relies on the results we obtained in \cite{CSS} 
on the existence and boundedness of imaginary powers of cone differential 
operators which -- according to a theorem of Dore and Venni \cite{DoVe} --
implies the maximal $L_{q}$-regularity. 
As a specific example we can treat the case where 
$A(u)=-a(t^c u)\Delta$ is the Laplace-Beltrami 
operator for a conic manifold of dimension greater than four, multiplied by
a  positive $\calC^\infty$-function $a$ depending on $t^c u$.
Here $t$ is a smooth function on $\bz$, 
which is strictly positive and extends the above coordinate $t$; 
$c$ is a positive constant.
The nonlinearity at the right hand side can
be taken to be a linear combination of functions of
the form $f(u)=|u|^\alpha$, $\alpha \ge 1$, or 
$f(u)=u^\alpha$, $\alpha \in \nz$.

The analysis on conic manifolds shows many interesting features.
One basic problem concerns the domain of the operators. On a closed 
manifold, an elliptic differential operator defined on all smooth 
functions has a single closed extension in $L_{p}$; its domain is the 
corresponding Sobolev space, which depends only on the order of the 
operator. For cone differential operators the situation is quite
different. They naturally
act on scales of 
weighted $L_{p}$-Sobolev spaces which coincide with the usual ones in 
the interior and are characterized by a weight function of the type 
$t^\gamma$, $\gamma \in \rz$, close to the boundary. 
For an elliptic operator, 
defined a priori on $\cicomp(\intb)$, there are, in general, many
different closed extensions,
parametrized by the subspaces of a finite-dimensional space of singular 
functions. 
They depend on the form of $A$ near $t=0$, as we shall see in 
Section \ref{semi3}. 
If one tries to employ maximal regularity techniques, the choice of the
domain therefore is of crucial importance. 

It is our intention to make the paper readable also for 
non-specialists in singular calculus. We shall highlight the specific 
difficulties of the subject and study many examples.

{\bf Acknowledgment:} We thank M.\ Hieber for helpful discussions  and 
his comments on an earlier version of the paper which lead to an improvement of the results.

\section{Differential operators on conic manifolds}\label{semi2}
In this section we recall some basic notions on cone differential operators 
and weighted Sobolev spaces. We summarize how ellipticity of such operators 
is described in terms of the symbolic structure, and how it is connected to 
the Fredholm property of the associated mapping between the Sobolev spaces.

\subsection{Operators of Fuchs type}\label{semi2.1}

Let $\bz$ be a compact manifold with boundary $X=\partial \bz$.
Fix once and for all a splitting of 
coordinates $(t,x)\in{[0,1[}\times X$ near the boundary $X$ of $\bz$. 
A {\em cone differential operator} or {\em Fuchs type operator} on $\bz$ 
is a differential operator -- or also a system of differential operators --
with
smooth coefficients on $\intb$ which near the 
boundary has the form 
  \begin{equation}\label{fuchs}
    A=t^{-\mu}\smsum_{j=0}^\mu a_{j}(t)(-t\partial_{t})^j,\qquad
    a_{j}\in\ci([0,1[,\text{\rm Diff}^{\,\mu-j}(X)),
  \end{equation}
where $\mu$ is the order of $A$. 
Let us stress 
the three main features of a Fuchs type operator: the singular 
factor $t^{-\mu}$ determined by the order of $A$, the smoothness of the 
coefficients $a_{j}$ up to $t=0$, and the totally characteristic 
derivatives in $t$-direction. Without difficulty, one could also 
treat a singular factor $t^{-\nu}$ for some real $\nu>0$.

Besides the usual homogeneous principal symbol 
  \begin{equation}\label{principal}
      \sigma^\mu_{\psi}(A)\in\ci((T^{*}\intb)\setminus 0)
  \end{equation}
taking values in bundle homomorphisms, we
associate with a cone differential operator $A$ two further symbolic 
levels: The {\em rescaled symbol} 
  \begin{equation}\label{rescaled}
      \tilde{\sigma}^\mu_{\psi}(A)\in\ci((T^{*}X\times\rz)\setminus 0)
  \end{equation}
is defined, in local terms, by 
  $$\tilde{\sigma}^\mu_{\psi}(A)(x,\xi,\tau)=
    \smsum_{j=0}^\mu\sigma^\mu_{\psi}(a_{j})(0,x,\xi)(-i\tau)^j.$$
The {\em conormal symbol} 
  \begin{equation}\label{conormal}
      \sigma^\mu_{M}(A)\in\calA(\cz,{\rm Diff}^\mu(X))\subset
      \calA(\cz,\calL(H^s_{p}(X),H^{s-\mu}_{p}(X)))
  \end{equation}
is an entire function taking values in (systems of)
differential operators on the 
boundary $X$. It is given by 
  $$\sigma^\mu_{M}(A)(z)=\smsum_{j=0}^\mu a_{j}(0)z^j.$$
Ellipticity of $A$ shall be described in terms of the invertibility 
of the symbols \eqref{principal}, \eqref{rescaled}, and 
\eqref{conormal}. 
As the case of systems of operators does not present additional 
analytical difficulties, we shall not stress this point in the text, below.
\begin{example}\label{laplace1}
    Let $g(t)$ be a family of smooth metrics on $X$ that depends 
    smoothly on a parameter $t\in[0,1[$. Equip $\intb$ with 
    a metric that coincides with $dt^2+t^2g(t)$ near $t=0$. Near the boundary, 
    the associated Laplace-Beltrami operator $\Delta$ is given by 
      $$t^{-2}\left\{(t\partial_{t})^{2}+
        (n-1+t(\log G)'(t))t\partial_{t}+\Delta_{X}(t)\right\},$$
    where $G=\det(g_{ij})^{1/2}$ and $\Delta_{X}(t)$ is the Laplacian 
    on $X$ with respect to the metric $g(t)$. Thus $\Delta$ is a second order 
    Fuchs type operator on $\bz$ with rescaled symbol 
      $$\tilde\sigma^2_{\psi}(\Delta)(x,\tau,\xi)=-\tau^{2}-|\xi|^{2},$$
    where $|\xi|$ refers to the metric $g(0)$. Its conormal symbol is
      $$\sigma^2_{M}(\Delta)(z)=z^{2}-(n-1)z+\Delta_{X}(0).$$
\end{example}
\subsection{Weighted cone Sobolev spaces}\label{semi2.2}
The intention to find a class of spaces on which Fuchs type operators are 
naturally continuous leads to the definition of the following scale of weighted 
Sobolev spaces on the interior of $\bz$: 
\begin{definition}\label{sobolev}
    Let $s\in\nz_{0}$, $\gamma\in\rz$, and $1<p<\infty$. Then $\hsgpb$ 
    denotes the space of all distributions 
    $u\in H^s_{p,\text{\rm loc}}(\intb)$ such that 
      $$t^{\frac{n+1}{2}-\gamma}(t\partial_{t})^k\partial_{x}^\alpha
        (\omega u)(t,x)\;\in\;L_{p}([0,1[\times 
        X,\mbox{$\frac{dt}{t}dx$})\qquad
	\forall\; k+|\alpha|\le s$$
    for some cut-off function $\omega$ (the particular choice of $\omega$ is irrelevant). 
\end{definition}
Recall that a cut-off function is a function $\omega\in\cicomp([0,1[)$ such that 
$\omega\equiv1$ near $t=0$. The index $s$ indicates the smoothness of functions, 
while the weight index $\gamma$ measures the flatness or rate of vanishing near 
the boundary. We shall extend this definition to arbitrary $s\in\rz$ in the sequel. 
A Fuchs type operator $A$ as in \eqref{fuchs} clearly induces continuous mappings 
 $$A:\hsgpb\to\calH^{s-\mu,\gamma-\mu}_{p}(\bz)$$ 
for any $\gamma$, $s$, and $p$. Apart from a certain normalization, the particular 
choice of the weight factor $t^{\frac{n+1}{2}}$ and the measure $\frac{dt}{t}dx$ 
ensures that essential properties of $A$, like Fredholm property and invertibility, 
are independent of $s\in\rz$ and $1<p<\infty$, cf.\ Theorem \ref{fredholm}. 
\begin{example}\label{lpspace}
    If we identify $\rz^{1+n}$ with $\rz_{+}\times S^n$ via polar coordinates 
    (i.e.\ we regard 0 as a conic singularity), a function $u$ 
    belongs to $L_{p}(\rz^{1+n})$ if and only if 
    $t^{\frac{n+1}{p}}u(t,x)\in L_{p}(\rz_{+}\times 
    S^n,\frac{dt}{t}dx)$. This suggests to regard 
       \begin{equation}\label{eq:lpspace}
	   L_{p}(B):=\calH^{0,\gamma_{p}}_{p}(\bz),\qquad
	   \gamma_{p}=(n+1)\left(\mbox{$\frac{1}{2}-\frac{1}{p}$}\right),
       \end{equation}
    as the natural $L_{p}$-spaces on the conic manifold $B$. 
\end{example}
Writing the standard Sobolev spaces $H^s_{p}(\rz^{1+n})$ in polar 
coordinates leads to more complicated spaces (so-called subspaces with 
asymptotics), cf.\ \cite{Daug}, Appendix A, \cite{Schu2}, Theorem 
1.1.22. We shall illustrate this later on, see \eqref{hpspace} in Example 
\ref{laplacian5}. 

There are various ways of extending the definition of cone Sobolev 
spaces to real smoothness parameters $s$, for example by interpolation 
and duality. For later purposes we want to sketch a definition based on 
the use of local coordinates. To this end let 
    \begin{equation}\label{local}
        H^{s,\gamma}_{p}(\rz^{1+n}_{(t,x)})=
	\left\{u\in\calD'(\rz^{1+n})\st
	e^{-\gamma\spk{t}}u(t,x)\in H^{s}_{p}(\rz^{1+n})\right\}
    \end{equation}
with the canonically induced norm. As usual, $\spk{t}=(1+t^2)^{1/2}$. Moreover 
let $S_{\gamma}:\calD'(\rz^{1+n})\to\calD'(\rz^{1+n})$ be defined by 
    \begin{equation}\label{sgamma}
	S_{\gamma}:\cicomp(\rz^{1+n})\to\cicomp(\rz^{1+n}),\quad
	u(t,x)\mapsto e^{(\frac{n+1}{2}-\gamma)t}u(e^{-t},x).
    \end{equation}
Let $\kappa_{j}:U_{j}\to\rz^n$, $j=1,\ldots,N$, and 
$\chi_{j}:V_{j}\to\rz^{1+n}$, $j=1,\ldots,M$, provide coverings by 
coordinate charts of $X$ and $\bz$, respectively, and 
$\{\varphi_{j}\}$, $\{\psi_{j}\}$ be corresponding subordinate 
partitions of unity. Then $\hsgpb$ is the space of all distributions 
such that 
    \begin{equation}\label{norm}
	\begin{array}{lcl}
	    \|u\|_{\hsgpb} \! & \! = \! & \!
	    \displaystyle
	    \smsum_{j=1}^{N}\|S_{0}(1\times\kappa_{j})_{*}
	              (\omega\varphi_{j}u)\|_{H^{s,\gamma}_{p}(\rz^{1+n})}+ 
	    \smsum_{j=1}^{M}\|\chi_{j*}
	              ((1-\omega)\varphi_{j}u)\|_{H^{s}_{p}(\rz^{1+n})}
	\\
		      \! & \! = \! & \!
		      \displaystyle
	    \smsum_{j=1}^{N}\|S_{\gamma}(1\times\kappa_{j})_{*}
	              (\omega\varphi_{j}u)\|_{H^{s}_{p}(\rz^{1+n})}+ 
	    \smsum_{j=1}^{M}\|\chi_{j*}
	              ((1-\omega)\varphi_{j}u)\|_{H^{s}_{p}(\rz^{1+n})}
        \end{array}    
    \end{equation}
is defined and finite. Here, $\omega\in\cicomp([0,1[)$ is a cut-off 
function and $*$ refers to the push-forward of distributions. 
Up to equivalence of norms, this construction is independent of the choice of $\kappa_j$ and $\chi_j$. 
\subsection{Ellipticity of cone differential operators}\label{semi2.3}
Each Fuchs type operator $A$ of the form \eqref{fuchs} induces 
continuous actions $A:\hsgpb\to\calH^{s-\mu,\gamma-\mu}_{p}(\bz)$ for 
any $\gamma,s\in\rz$ and $1<p<\infty$. We next address the question when it is a 
Fredholm operator. 

A cone differential operator $A$ is called {\em elliptic} with respect 
to the weight $\gamma\in\rz$ if the following conditions are 
satisfied: 
    \begin{itemize}
	\item[(1)] Both the homogeneous principal symbol 
	     $\sigma_{\psi}^\mu(A)$ and the rescaled symbol 
	     $\tilde\sigma_{\psi}^\mu(A)$ are invertible, 
	\item[(2)] the conormal symbol is invertible on the line $\re 
	     z=\frac{n+1}{2}-\gamma$, i.e. 
	       $$\sigma_{M}^\mu(A)(z):H^s_{p}(X)\xrightarrow{\;\simeq\;} 
	         H^{s-\mu}_{p}(X)\qquad
	         \forall\;\re z=\mbox{$\frac{n+1}{2}-\gamma$}.$$
    \end{itemize} 
Due to the spectral invariance of pseudodifferential operators on 
closed manifolds, condition (2) is independent of the 
choice of $s$ and $p$. 

Under condition (1) the conormal symbol 
$\sigma_{M}^\mu(A)(z)\in\calA(\cz,L^\mu(X)))$ 
is meromorphically invertible with only finitely many poles in each 
vertical strip $|\re z|<k$, $k\in\nz$, cf.\ \cite{Schu2}, Theorem 2.4.20. 
Condition (2) is imposed in order to ensure that none of these poles lies on the line 
$\re z=\frac{n+1}{2}-\gamma$. 

Together, (1) and (2) allow the construction of a parametrix, cf.\ 
\cite{Schu2}, \cite{Melr}.  

The following theorem was shown in \cite{ScSe}: 
\begin{theorem}\label{fredholm}
    Let $A$ be a cone differential operator. Then the operator 
    $A\!:\!\hsgpb\to\calH^{s-\mu,\gamma-\mu}_{p}(\bz)$ is Fredholm if and 
    only if $A$ is elliptic with respect to the weight $\gamma$. The Fredholm property 
    as well as the index are independent of $s$ and $p$.  
\end{theorem}
\begin{example}\label{laplace2}
    Let $\Delta$ be the Laplacian on $\bz$ as described in Example
    \text{\rm \ref{laplace1}}. If $0=\lambda_{0}>\lambda_{1}>\ldots$ 
    are the eigenvalues of $\Delta_{X}(0)$, then 
    $\sigma_{M}^{2}(\Delta)(z)=z^2-(n-1)z+\Delta_{X}(0)$ is not 
    bijective if and only if 
      $$z\in\left\{\mbox{$\frac{n-1}{2}\pm
        \big(\frac{(n-1)^{2}}{4}-\lambda_{j}\big)^{1/2}$}
	\st j\in\nz_{0}\right\}.$$
    Accordingly, $\Delta$ is elliptic with respect to all $\gamma$ not 
    belonging to this set. For later purpose let us point out that in 
    any case $\sigma_{M}^{2}(\Delta)$ is invertible in the strip 
    $0<\re z<n-1$. 
\end{example}

\section{Closed extensions of cone differential operators}\label{semi3}

We consider $A$ as an unbounded operator in $\calH^{0,\gamma}_{p}(\bz)$, 
\begin{equation}\label{unbounded}
    A:\cicomp(\intb)\subset\calH^{0,\gamma}_{p}(\bz)\longrightarrow
    \calH^{0,\gamma}_{p}(\bz)
\end{equation} 
and shall investigate its closed extensions. The material in Proposition \ref{dmin} 
through Corollary \ref{unique} goes back to Lesch's work \cite{Lesc} for the 
case $p=2$ and we omit proofs. We shall assume that $A$ is elliptic in the interior, 
i.e.\ satisfies the ellipticity condition (1) of Section \ref{semi2.3}. 
    
In contrast to elliptic pseudodifferential operators on closed manifolds, a cone 
differential operator $A$ has in general infinitely many closed extensions. There 
are two natural extensions - the minimal and maximal extension 
$\amin=A^{\gamma,p}_{\min}$ and $\amax=A^{\gamma,p}_{\max}$. The minimal extension 
is the closure of the operator $A$ in \eqref{unbounded}, hence   
  $$\calD(\amin)=\left\{u\in\calH^{0,\gamma}_{p}(\bz)\st 
         \exists\,(u_{n})\subset\cicomp(\intb):\ 
         u_{n}\to u\text{ and }Au_{n}\to v=:Au\text{ in }
	 \calH^{0,\gamma}_{p}(\bz)\right\};$$
the maximal extension is given by the action of $A$ on the domain 
  $$\calD(\amax)=\left\{u\in\calH^{0,\gamma}_{p}(\bz)\st Au\in
         \calH^{0,\gamma}_{p}(\bz)\right\}.$$
These two special cases are the key to understanding the general situation. 
\begin{proposition}\label{dmin}
    The domain of the closure of $A$ is given by         
      $$\calD(\amin)=
        \Big\{u\in\mathop{\mbox{\Large$\cap$}}_{\varepsilon>0}
        \calH^{\mu,\gamma+\mu-\varepsilon}_{p}(\bz)\st 
        Au\in\calH^{0,\gamma}_{p}(\bz)\Big\}.$$
    In particular, 
        $$\calH^{\mu,\gamma+\mu}_{p}(\bz)\hookrightarrow\calD(\amin)
	\hookrightarrow
	  \calH^{\mu,\gamma+\mu-\varepsilon}_{p}(\bz)
	  \qquad\forall\;\varepsilon>0.$$ 
    If $A$ additionally satisfies condition \text{\rm (2)} of Section 
    \text{\rm{\ref{semi2.3}}} with respect to the weight 
$\gamma+\mu$, then 
    (topologically)
        $$\calD(\amin)=\calH^{\mu,\gamma+\mu}_{p}(\bz).$$
\end{proposition} 
If the coefficients $a_{j}$ in \eqref{fuchs} are independent of $t$ for $t$ close 
to $0$, this result follows from a simpler version of the above mentioned parametrix 
construction. The general case can be treated by means of perturbation theory, 
since for 
    $$A^{(0)}=\omega\, 
      t^{-\mu}\smsum_{j=0}^{\mu}a_{j}(0)(-t\partial_{t})^j\,\omega_{0}+
      (1-\omega)\,A\,(1-\omega_{1})$$
the cut-off functions $\omega,\omega_{j}\in\cicomp([0,1[)$ can be 
chosen in such a way that $A-A^{(0)}$ is $A^{(0)}$-bounded with 
$A^{(0)}$-bound less than 1, hence $\calD(\amin)=\calD(A^{(0)}_{\min})$, 
cf.\ \cite{Kato}, Theorem 1.1 on page 190.  

Assuming merely the interior ellipticity of $A$, one obtains:
\begin{proposition}
    $\amin$ is a Fredholm operator.
\end{proposition} 
Let us now turn to the description of the maximal extension of $A$. 
As mentioned in Section \ref{semi2.3}, the conormal symbol of $A$ is 
meromorphically invertible. Let $p_{1},\ldots,p_{N}$ denote the 
finitely many poles in the strip $\frac{n+1}{2}-\gamma-\mu<\re z<
\frac{n+1}{2}-\gamma$. Near each $p_{j}$ we write 
    \begin{equation}\label{laurent}	      
	\sigma_{M}^\mu(A)(z)^{-1}\equiv\smsum_{k=0}^{m_{j}}R_{jk}(z-p_{j})^{-k-1} 
    \end{equation} 
modulo a function holomorphic near $p_j$. It can be shown that 
the Laurent coefficients $R_{jk}$ belong to $L^{-\infty}(X)$ and have 
finite 
dimensional range. Define 
    $$G_{A}=G_{A}^{\gamma,p}=\text{\rm diag}(G_{1},\ldots,G_{N}): 
      \mathop{\mbox{\Large$\oplus$}}_{j=1}^{N}\ci(X)^{m_{j}+1}
      \longrightarrow
      \mathop{\mbox{\Large$\oplus$}}_{j=1}^{N}\ci(X)^{m_{j}+1}$$
by the left upper triangular matrices
    $$G_{j}=(g^j_{ik})_{0\le i,k\le m_{j}}\quad\text{ with }\quad 
      g^j_{ik}=\begin{cases}
                  R_{j,i+k} & \text{if }i+k\le m_{j}\\
		  0 & \text{else}
             \end{cases}.$$
$G_{A}$ is a finite rank operator. 
\begin{proposition}\label{dmax}
    There exists a finite dimensional vector space 
    $\calE_{A}=\calE_{A}^\gamma\subset\calH^{\infty,\gamma}_{p}(\bz)$ 
with 
        $$\calD(\amax)=\calD(\amin)\oplus\calE_{A},\qquad
	  \dim\calE_{A}=\text{\rm rank}\,G_{A},$$
	  as a topologically direct sum.
    The space $\calE_{A}$ does not depend on $1<p<\infty$. 
\end{proposition} 
\begin{corollary}\label{unique}
    Let 
    $A:\cicomp(\intb)\subset\calH^{0,\gamma}_{p}(\bz)\to\calH^{0,\gamma}_{p}(\bz)$ 
    be given. Then: 
        \begin{itemize}
            \item[a)] Any closed extension of $A$ is given by the action 
                of $A$ on a domain $\calD(\amin)\oplus V$ with a subspace 
                $V$ of $\calE_{A}$. 
            \item[b)] $A$ has a unique closed extension 
                $A^{\gamma,p}_{\min}=A^{\gamma,p}_{\max}$ if and only if the conormal 
                symbol $\sigma_{M}^\mu(A)(z)$ is invertible for all $z$ with 
                $\frac{n+1}{2}-\gamma-\mu<\re z<\frac{n+1}{2}-\gamma$. 
        \end{itemize}
\end{corollary}
\begin{example}\label{laplace3}
    If $\Delta$ is the Laplacian on $\bz$ as introduced in Example 
    \text{\rm\ref{laplace1}}, we saw in Example 
\text{\rm\ref{laplace2}} 
    that the conormal 
    symbol $\sigma_{M}^2(\Delta)(z)$ is invertible for 
    all $z$ with $0<\re z<n-1$ but not for $z=0$ and $z=n-1$. Hence
        $$\Delta:\cicomp(\intb)\subset\calH^{0,\gamma_{p}}_{q}(\bz)
	  \longrightarrow\calH^{0,\gamma_{p}}_{q}(\bz),\qquad 
	  1<p,q<\infty,$$
    cf.\ \eqref{eq:lpspace}, has a unique closed extension if and 
    only if $\frac{n+1}{2}-\gamma_{p}-2\ge0$ and 
    $\frac{n+1}{2}-\gamma_{p}\le n-1$. These conditions are satisfied 
    if and only if 
        \begin{equation}\label{pbound}
	    2\max(p,p')\le n+1,
	\end{equation}
    where, as usual, $p'$ denotes the number dual to $p$, i.e.\ 
    $\frac{1}{p}+\frac{1}{p'}=1$. 
\end{example}
To describe the space $\calE_{A}$ from Proposition \ref{dmax}, let us 
assume for simplicity that the coefficients $a_{j}$ in 
\eqref{fuchs} are independent of $t$ for $t$ close to $0$. 
If the inverted conormal symbol of $A$ is as in \eqref{laurent}, then 
    \begin{equation}\label{asymptotics}
	   \calE_{A}=\Big\{
	                \omega\,\smsum_{j=1}^{N}\smsum_{l=0}^{m_{j}}
		        \zeta_{jl}(u)t^{-p_{j}}\log^lt\st
		        u\in\cicomp(]0,1[\,\times X)
	          \Big\}
    \end{equation}
with an arbitrary (fixed) cut-off function $\omega\in\cicomp([0,1[)$ 
and the linear finite rank mappings $\zeta_{jl}:\cicomp(]0,1[\,\times 
X)\to\ci(X)$ being defined by 
$$\zeta_{jl}(u)=\frac{(-1)^l}{l!}\smsum_{k=l}^{m_j}\frac{1}{(k-l)!}      
  R_{jk}\frac{\partial^{k-l}(\calM u)}{\partial z^{k-l}}(p_j+\mu);$$
here, $\calM u=\calM_{t\to z}u\in\calA(\cz,\ci(X))$ denotes the Mellin 
transform of $u$. In case the coefficients depend on $t$, one can show that 
$\calE_A\subset V$, where $V$ is a finite-dimensional space of singular functions 
which is similar to the right hand side in \eqref{asymptotics}.  
\begin{example}\label{laplace4}
    Let us reconsider the Laplacian introduced in Example 
    \text{\rm\ref{laplace1}} specializing to $\dim\bz=2$ with $X=S^1=\rz/2\pi\gz$ 
and metric 
    $dt^{2}+t^{2}d\vartheta^{2}$ on $[0,1[\times S^1$, where 
$d\vartheta^{2}$ is the standard 
    metric on $S^{1}$. Then the conormal symbol 
    $\sigma_{M}^{2}(\Delta)(z)=z^{2}+\partial_{\vartheta}^{2}$ has the 
    non-bijectivity points $p_{j}=j\in\gz$. Passing to Fourier 
    series, outside $\gz$ the inverse is given by 
        $$(z^{2}+\partial_{\vartheta}^{2})^{-1}:\quad
          \smsum_{k}c_{k}e^{ik\vartheta}\mapsto
          \smsum_{k}\mbox{$\frac{c_{k}}{z^2-k^{2}}$}e^{ik\vartheta}.$$
    For fixed $j$, only the terms coming from $k=\pm j$ are not 
    holomorphic near $j$. This shows that
        $$(z^{2}+\partial_{\vartheta}^{2})^{-1}\equiv R_{01}\,z^{-2}$$
    near $z=0$ modulo holomorphic functions, while, near $z=j\not=0$,
        $$(z^{2}+\partial_{\vartheta}^{2})^{-1}\equiv 
R_{j0}\,(z-j)^{-1},$$
    where the Laurent coefficients are given by 
        $$R_{01}f(\vartheta)=\hat{f}_{0},\qquad
	        R_{j0}f(\vartheta)=\mbox{$\frac{1}{2j}$}\,\hat{f}_{-j}\,e^{-ij\vartheta}+
	        \mbox{$\frac{1}{2j}$}\,\hat{f}_{j}\,e^{ij\vartheta}$$
    ($\hat{f}_{k}$ denoting the $k$-th Fourier coefficient of $f$). 
    Thus we obtain 
        $$\dim\calD(\Delta_{\max}^{\gamma,p})/
	  \calD(\Delta_{\min}^{\gamma,p})=
	      \begin{cases}
		  2 & \text{if }\gamma\in\gz \\
		  4 & \text{else}
	      \end{cases}.$$ 
    In the particular case $\gamma_{p}=1-\frac{2}{p}$, cf.\ 
    \text{\rm \eqref{eq:lpspace}}, the domain of the maximal 
extension 
    is given by 
        $$\calD(\Delta_{\max}^{0,q})=\calD(\Delta_{\min}^{0,q})\oplus 
	       \omega\,\text{\rm span}(1,\log t),\qquad p=2,$$
    for any $1<q<\infty$, and for $p\not=2$ by 
        $$\calD(\Delta_{\max}^{\gamma_{p},q})=
	  \calH_{q}^{2,\gamma_{p}+2}(\bz)\oplus\omega\,
	      \begin{cases}
		  \text{\rm span}(1,\log t,e^{i\vartheta}t,e^{-i\vartheta}t) &
		      p>2 \\
		  \text{\rm span}(1,\log t,e^{i\vartheta}t^{-1},
		                  e^{-i\vartheta}t^{-1}) &
		      p<2
	      \end{cases}.$$ 
\end{example}

\section{Bounded imaginary powers}\label{semi4}

The boundedness of purely imaginary powers $A^{iy}$, $y\in\rz$, of an 
operator $A:\calD(A)\subset Y\to Y$ is closely related to the unique 
solvability of the parabolic equation 
    \begin{equation}\label{parabolic}
	\begin{cases}
	    \dot{u}+Au=f & \mbox{on $]0,T[$}
	    \\
	    u(0)=u_{0}.
	\end{cases}
    \end{equation}
In \cite{DoVe} Dore and Venni proved the following 
theorem: 
\begin{theorem}\label{dorevenni}
    Let $A:\calD(A)\subset Y\to Y$ be a closed densely defined and 
    positive operator in a \text{\rm UMD}-space $Y$. If the imaginary 
    powers of $A$ exist and satisfy the estimate 
        \begin{equation}\label{bip}
	    \|A^{iy}\|_{\calL(Y)}\le c\,e^{\theta|y|}\qquad \forall\;y\in\rz
	\end{equation}
    for some $0<\theta<\frac{\pi}{2}$, then the initial value 
    problem \eqref{parabolic} with $u_{0}=0$ has, for any $f\in 
L_{r}([0,T],Y)$,  $1<r<\infty$, a unique solution 
        $$u\in W^{1}_{r}([0,T],Y)\;\cap\;L_{r}([0,T],\calD(A)).$$
    Moreover, $u$, $\dot{u}$ and $Au$ depend continuously on $f$. 
\end{theorem} 
Positivity of a linear operator here means that the resolvent set $\varrho(A)$ contains all non-negative reals, and $\|(A+\lambda)^{-1}\|_{\calL(Y)}=O(\lambda^{-1})$ for $\lambda\ge0$. 
In applications, the assumption on $Y$ to be a UMD-space is not very 
restrictive. For example, $L_{p}(\Omega,d\mu)$, $1<p<\infty$, is a 
UMD-space for any $\sigma$-finite measure space $(\Omega,\mu)$, cf.\ 
\cite[Theorem 4.5.2]{Aman}. This is then also true for the cone Sobolev spaces 
$\hsgpb$, since $\calH^{0,\gamma}_{p}(\bz)$ is a weighted 
$L_{p}$-space on $\bz$ and, by the existence of order reductions, 
$\hsgpb$ is isomorphic to $\calH^{0,\gamma}_{p}(\bz)$ for any 
$s\in\rz$. 

The key assumption of Theorem \ref{dorevenni} is the existence of the 
imaginary powers together with the estimate \eqref{bip}. In 
\cite{CSS} 
we gave criteria, when this assumption holds true for the minimal or 
maximal extension of a cone differential operator $A$. To describe 
these criteria let us recall the notion 
of the {\em model cone operator} $\widehat{A}$ associated with $A$. 
The idea is to freeze the coefficients of $A$ at $t=0$ so that we obtain an 
operator that lives on the infinite cone over $X$. On this cone we have a 
natural choice of Sobolev spaces together with a weight function at the origin. 
For practical reasons, we work on the cylinder $X^\wedge:=\rz_{+}\times X$ 
with $(t,x)$-coordinates. Then 
    \begin{equation}\label{ahat}
	\widehat{A}=t^{-\mu}\smsum_{j=0}^\mu a_{j}(0)(t\partial_{t})^{j}\;:\;
	\cicomp(X^\wedge)\subset\calK^{0,\gamma}_{p}(X^\wedge)\longrightarrow
	\calK^{0,\gamma}_{p}(X^\wedge), 
    \end{equation}
for $A$ as in \eqref{fuchs} and the scale of Sobolev spaces is defined as follows: 
\begin{definition}\label{ksg}
    $\calK^{s,\gamma}_{p}(X^\wedge)$ consists of all distributions 
    $u\in H^{s}_{p,\text{\rm loc}}(X^\wedge)$ such that for some 
    cut-off function $\omega\in\cicomp([0,1[)$ 
        \begin{itemize}
	    \item[i)] $\omega u\in\hsgpb$,
	    \item[ii)] if $\kappa:U\to\rz^{n}$ is a coordinate chart of $X$ 
	        and $\chi(t,x)=(t,t\kappa(x))$, then 
	        $\chi_{*}[(1-\omega)\phi u]\in H^{s}_{p}(\rz^{1+n})$ for 
	        any $\phi\in\cicomp(U)$.
	\end{itemize}
\end{definition}
For $p=2$, the spaces $\calK^{s,\gamma}_{p}(X^\wedge)$ were introduced by Schulze, see \cite{Schu2}. 

If $A$ satisfies ellipticity condition (1) of Section \ref{semi2.3}, 
the extensions of $\widehat{A}$ can be described quite similar to the 
extensions of $A$. In particular, 
$$\calD(\widehat{A}_{\max})=\calD(\widehat{A}_{\min})\oplus\calE_{A}$$
with $\calE_{A}$ from \eqref{asymptotics} (now considered as a function 
space on $X^\wedge$). In case $A$ also satisfies 
condition (2) of Section \ref{semi2.3} with respect to the weight 
$\gamma+\mu$,
$$\calD(\widehat{A}_{\min})=\calK^{\mu,\gamma+\mu}_{p}(X^\wedge).$$
\begin{example}\label{laplacian5}
    Let $\Delta$ be the Laplacian on $\bz$ from Example 
    \text{\rm\ref{laplace4}}. Then $\widehat{\Delta}$ is given by 
    $t^{-2}((t\partial_{t})^2+\partial_{\vartheta}^2)$ acting on 
    $\calK^{0,\gamma}_{p}(S^{1\wedge})$. In the special case 
    $\gamma=\gamma_{p}$ we have 
    $\calK^{0,\gamma_{p}}_{p}(S^{1\wedge})=L_{p}(\rz^2)$ via polar 
    coordinates, cf.\ Example \text{\rm\ref{lpspace}}. Hence the 
    associated model cone operator is just the standard 
    Laplacian, i.e. 
        $$\widehat{\Delta}=\widehat{\Delta}^p=\Delta_{\rz^2}\;:\;
	  \cicomp(\rz^2\setminus\{0\})\subset
	  L_{p}(\rz^2)\longrightarrow L_{p}(\rz^2).$$
    As in Example \text{\rm\ref{laplace4}} we get 
      $$\calD(\widehat{\Delta}_{\max}^p)=\calD(\widehat{\Delta}_{\min}^p)
	      \oplus\omega\,
        \begin{cases}
	        \text{\rm span}(1,\log t,e^{i\vartheta}t,e^{-i\vartheta}t) & p>2\\
	        \text{\rm span}(1,\log t) & p=2\\
	        \text{\rm span}(1,\log t,e^{i\vartheta}t^{-1},e^{-i\vartheta}t^{-1}) & p<2
	      \end{cases}$$
    with $\calD(\widehat{\Delta}_{\min}^p)=
    \calK^{2,\gamma_{p}+2}_{p}(S^{1\wedge})$ in case $p\not=2$. A 
    particular closed extension of $\widehat{\Delta}^p$ is given by 
    the action of $\Delta_{\rz^2}$ on $H^{2}_{p}(\rz^2)$. We then 
    obtain 
        \begin{equation}\label{hpspace}            
H^{2}_{p}(\rz^2)=\calD(\widehat{\Delta}_{\min}^p)\oplus\omega\,
	    \begin{cases}
	        \text{\rm span}(1,e^{i\vartheta}t,e^{-i\vartheta}t) & p>2\\
	        \text{\rm span}(1) & p\le 2
            \end{cases}.
        \end{equation}
    In fact, we know that 
    $H^{2}_{p}(\rz^2)\subset\calD(\widehat{\Delta}_{\max}^p)$. Moreover, the Sobolev 
    embedding theorem ensures that 
    $H^{2}_{p}(\rz^2)$ is a subspace of $\calC(\rz^2)$. Thus, in case $p>2$, the 
    function $\log t$ cannot belong to $H^{2}_{p}(\rz^2)$. For the same reason, we can 
    exclude $\text{\rm span}(\log t,e^{i\vartheta}t^{-1},e^{-i\vartheta}t^{-1})$ in 
    case $p\le 2$. Finally, in cartesian coordinates, 
    $e^{\pm i\vartheta}t=x\pm iy$ is a smooth function on $\rz^2$ and thus 
    $\omega e^{\pm i\vartheta}t$ belongs to $H^{2}_{p}(\rz^2)$. 
\end{example}
\begin{theorem}\label{criteria}
    Let the cone differential operator $A$ satisfy 
        \begin{itemize}
            \item[\text{\rm(E)}] 
	      both $\sigma_{\psi}^\mu(A)$ and $\tilde{\sigma}_{\psi}^\mu(A)$ 
	      have no spectrum in $\Lambda_{\theta}$, 
	\end{itemize}
    where $\Lambda_{\theta}=
    \{z\in\cz\st|\text{\rm arg}z|\ge\theta\}\cup\{0\}$ is a closed 
sector and 
    $0\le\theta<\pi$.  
        \begin{itemize}
	    \item[a)] If $A$ satisfies condition \text{\rm(2)} of Section 
	        \text{\rm\ref{semi2.3}} with respect to $\gamma+\mu$ and 
		\\[1ex]
                    \hspace*{3mm}
		    $\mathrm{(E_{min})}$
			$\widehat{A}_{\min}$ has no spectrum in 
                        $\Lambda_{\theta}\setminus\{0\}$ \\[1ex]
		then there exists a $\varrho\ge0$ such that 
		    $$\|(\amin+\varrho)^{it}\|_{\calL(\calH^{0,\gamma}_{p}(\bz))}
		      \le c\,e^{\theta|t|}\qquad\forall\;t\in\rz.$$
	    \item[b)] If $A$ satisfies condition \text{\rm(2)} of Section 
	        \text{\rm\ref{semi2.3}} and 
		\\[1ex]
                    \hspace*{3mm}
		    $\mathrm{(E_{max})}$
                          $\widehat{A}_{\max}$ has no spectrum in 
                          $\Lambda_{\theta}\setminus\{0\}$ \\[1ex]
		then there exists a $\varrho\ge0$ such that 
		    $$\|(\amax+\varrho)^{it}\|_{\calL(\calH^{0,\gamma}_{p}(\bz))}
		      \le c\,e^{\theta|t|}\qquad\forall\;t\in\rz.$$	        
	\end{itemize}
\end{theorem}
Note that Theorem \ref{dorevenni} then also holds true for $A_{\min}$ 
and $A_{\max}$, respectively, since equation \eqref{parabolic} is 
equivalent to $\dot{v} + (A+\varrho) v = g$ with $g(t) = e^{\varrho t} f(t)$ and 
$v(t) = e^{\varrho t} u(t)$.

\begin{example}\label{exlaplmaxreg}
    Let $\Delta = \Delta^{\gamma_{p},q}$ be as in Example
    \text{\rm\ref{laplace3}}. If $2\max(p,p^\prime) - 1 < n$, cf. 
    \text{\rm\eqref{pbound}}, then $-\Delta$ satisfies both 
conditions {\rm     (E)} 
and $\mathrm{(E_{min})}$ for any $0 < \theta < \frac{\pi}{2}$. 
In fact, condition {\rm (E)} is clearly fulfilled. It is more difficult to
check ${\mathrm
(\rm E_{\mathrm{min}})}$. Details can be found in  
Theorem \text{\rm 7.1} of \cite{CSS}. 
    Correspondingly, the heat equation $\partial_{t} u - \Delta u = 
    f$, $u(0) = 0$, has a unique 
    solution $u \in W^{1}_{r}([0,T],\calH^{0,\gamma_{p}}_{q}(\bz)) 
    \cap L_{r}([0,T],\calH^{2,\gamma_{p}+2}_{q}(\bz))$ for any $f \in 
    L_{r}([0,T],\calH^{0,\gamma_{p}}_{q}(\bz))$, $1<q<\infty$.
\end{example}

The idea of proving Theorem \ref{criteria} is to consider $\lambda - 
A$ as an element of a suitable parameter-dependent pseudodifferential 
calculus on $\bz$ (i.e., the cone algebra as introduced by Schulze 
\cite{Schu2}). This technique was also used by Gil \cite{Gil}. 

Conditions (E) and $(\mathrm{E_{\min}})$ then assure 
that $\lambda - A$ is an elliptic element in this calculus, and 
therefore we find a parametrix $R(\lambda)$ which in fact coincides 
with the resolvent $(\lambda - A)^{-1}$ for large $\lambda \in 
\Lambda_{\theta}$. This yields, cf. Proposition 4.7 of \cite{CSS}, 
that $\| (\lambda-A)^{-1} \|_{\calL(\calH^{0,\gamma}_{p}(\bz))} =
\| R(\lambda) \|_{\calL(\calH^{0,\gamma}_{p}(\bz))} = O 
(|\lambda|^{-1})$ for $|\lambda| \to +\infty$ and thus complex powers 
of $A$ can be defined by the Dunford integral
\begin{equation}\label{dunfordint}
    A^z = \frac{1}{2\pi i} \int_{\calC} \lambda^z (\lambda-A)^{-1} 
    d\lambda, \hspace{1cm} \re z < 0,
\end{equation}
where $\calC$ is an appropriate path that coincides with $\partial 
\Lambda_{\theta}$ away from $0$. Assuming that $(\lambda - A)^{-1}$ 
exists in the whole sector $\Lambda_{\theta}$ (which can be achieved 
replacing $A$ by $A+c$), the use of the microlocal structure of 
$(\lambda - A)^{-1} = R(\lambda)$ allows to show that
$\| A^z \|_{\calL(\calH^{0,\gamma}_{p}(\bz))} \le c_{p} e^{\theta 
|\im 
z|}$ for $| \re z |$ sufficiently small. This estimate then extends 
to the purely imaginary powers. The result for the maximal extension 
follows from the one for the closure by passing to the adjoint.

The proof of Theorem \ref{criteria} relies only
on the structure of the parametrix to $\lambda - A$. Thus, corresponding results
are true for others than the  minimal or maximal extension of $A$, as soon as 
one finds criteria that ensure the existence of such a parametrix. 

The following example shows that, for the two-dimensional Laplacian, it is 
neither the minimal nor the maximal extension which is most interesting.
Instead we show that there is an intermediate extension generating a
holomorphic semigroup. 

Let $A = - \Delta$ be the Laplacian as in Example \ref{laplace4}. Denote 
by $A_{p}$ the extension with domain
    \begin{equation}\label{dap}
	\calD(A_{p}) = \calD(\Delta^{\gamma_{p},p}_{\mathrm{min}})
	\oplus \omega\,	       \begin{cases}
	           \mathrm{span}(1,t e^{i\vartheta}, t e^{-i\vartheta}), 
		   & p  >  2
		   \\
	           \mathrm{span}(1), 
		   & p \le 2		   
	       \end{cases}
    \end{equation}
    cf. \eqref{hpspace}. For $p \le 2$, the functions $\omega
    e^{i\vartheta} t$ 
    and $\omega e^{-i\vartheta} t$ both are elements of 
    $\calD(\Delta^{\gamma_{p},p}_{\mathrm{min}})$. In 
    fact, this follows from Proposition \ref{dmin}, since $\omega e^{\pm 
    i\vartheta}t \in \calH^{\infty,\gamma_{p}+2-\varepsilon}_{p}(\bz)$ for 
    every $\varepsilon > 0$, and $\Delta(e^{\pm i\vartheta} t) = 0$ implies that 
    $\Delta(\omega e^{\pm i\vartheta} t) \in
    \calH^{\infty,\gamma_{p}}_{p}(\bz)$. For 
    $p<2$, we even have $\omega \, t \, u \in 
    \calH^{\infty,\gamma_{p}+2}_{p}(\bz)$ for every function $u$ which is smooth 
    up to the boundary of $\bz$; this is no longer true in case $p=2$.
    
    Let $0 < \theta < \frac{\pi}{2}$ be arbitrary. 
    We shall show that the resolvent $(\lambda-A_{p})^{-1}$ exists for 
    all but finitely many $\lambda \in \Lambda_{\theta}$ and satisfies 
    \[
    \| (\lambda-A_p)^{-1} \|_{\calL(\calH^{0,\gamma}_{p}(\bz))} =
    O(|\lambda|^{-1}) \text{ for $|\lambda| \to \infty$}.
    \]
\begin{corollary}
    $-\Delta$ fulfills neither condition 
    $\mathrm{(E_{\min})}$ of Theorem {\rm\ref{criteria}} nor 
    $\mathrm{(E_{\max})}$. Even more is true: both $\lambda + 
    \Delta^{\gamma_{p},p}_{\mathrm{min}}$ and $\lambda + 
    \Delta^{\gamma_{p},p}_{\mathrm{max}}$ are non-invertible for 
    all $\lambda \in \cz$.
\end{corollary}
\begin{proof}
    Due to the compact embedding $\calD(\Delta^{\gamma_{p},p}_{\mathrm{max}}) 
    \hookrightarrow \calH^{0,\gamma_{p}}_{p}(\bz)$, the spectrum of     $\Delta^{\gamma_{p},p}_{\mathrm{max}}$ 
    is either all of $\cz$ or a discrete set. In the second case we thus would
    find a point $\lambda$, which belongs simultaneously to the resolvent sets of $A_p$ and
    $\Delta^{\gamma_{p},p}_{\mathrm{max}}$. However, this cannot be true, since 
    $\calD(A_p)$ is a proper subspace of the domain of the maximal extension. The 
    argument for the minimal operator is analogous, since its domain is a
    proper subspace of $\calD(A_p)$. 
\end{proof}
To obtain the statement on the resolvent of $A_p$,  let 
    $0 < \varrho < 1$ be fixed and let $\Delta_{2\bz}$ denote the 
    Laplacian on $2\bz$ (which is the double of $\bz$ or any smooth 
    closed manifold containing $\bz$ as a submanifold) with respect to a 
    metric that coincides with the given metric on $\bz 
    \setminus ([0,\frac{\varrho}{2}]\times X)$. Then define     
    $R_{p}(\lambda):\calH^{0,\gamma_p}_p(\bz)\to\calD(A_p)$ by     %
    \begin{equation}\label{Res1}
	R_{p}(\lambda) = \omega (\lambda - \widehat{A}_{p})^{-1} \omega_{0} +
	(1-\omega) (\lambda - \Delta_{2\bz}^p)^{-1}(1-\omega_{1})
    \end{equation}
    where $\omega, \omega_{j} \in C^\infty_{0}([0,1[)$ are cut-off 
    functions satisfying $\omega \equiv \omega_{j} \equiv 1$ on 
    $[0,\varrho]$ and $\omega \omega_{1} = \omega_{1}$, $\omega_{0} 
    \omega = \omega$. Moreover, let $\Delta^p_{2\bz}$ be the 
unbounded 
    operator in $L_{p}(2\bz)$ acting like $\Delta_{2\bz}$ on the 
    domain $H^{2}_{p}(2\bz)$. Note that \eqref{Res1} makes sense for 
    sufficiently large $\lambda \in \Lambda_{\theta}$ and
    \[
    \omega (\lambda - \widehat{A}_{p})^{-1} \omega_{0} =
    \omega\, {\mathrm{op}}_{y}\Big( \frac{1}{\lambda-|\eta|^{2}} \Big) 
    \omega_{0}
    \qquad\forall\; \lambda \notin [0,\infty[.
    \]
    On the right-hand side, {\rm op} denotes the usual pseudodifferential action,  
    $y$ and $\eta$ are the variables and covariables for $\rz^2$, if we identify 
    $L_{p}(\rz^{2})$ with 
    $\calK^{0,\gamma_{p}}_{p}(S^{1\wedge})$ and 
    $\calD(\widehat{A}_{p})$ with $H^{2}_{p}(\rz^{2})$ via polar 
    coordinates, cf.\ \eqref{hpspace}.
    
    The next lemma states that choosing other cut-off functions in 
    \eqref{Res1} changes $R_{p}(\lambda)$ only modulo ``good'' 
    remainders.
    
    \begin{lemma}\label{cutoff}
	Let $\sigma, \sigma_{0}, \sigma_{1} \in C^\infty_{0}([0,1[)$ be 
	arbitrary cut-off functions with $\sigma \sigma_{1} = \sigma_{1}$, 
	$\sigma_{0}\sigma = \sigma$ and $\sigma \equiv \sigma_{j} \equiv 1$ 
	on $[0,\varrho]$. Then
	\[
	R_{p}(\lambda) = \sigma (\lambda - \widehat{A}_{p})^{-1} \sigma_{0}
	               + (1-\sigma) (\lambda - \Delta^p_{2\bz})^{-1} 
(1-\sigma_{1})
		       + G(\lambda)
	\]
	with a remainder $G(\lambda)$
	\[
	G(\lambda) \in \calS(\tilde{\Lambda}_{\theta}, 
	\calL(\calH^{0,\gamma_{p}}_{p}(\bz), \calD(A_{p})))
	\]
	where $\tilde{\Lambda}_{\theta} = \{ z \in \Lambda_{\theta} \; | \; 
|z| 	\ge 
c\}$, $0 < \theta < \frac{\pi}{2}$, with a sufficiently large constant
	$c>0$.
    \end{lemma}
    \begin{proof}
	Let us write $\ci(\bz)$ for the space of functions which are smooth up to the 	
	boundary of $\bz$. Since the scalar product 	
	$\skp{\cdot}{\cdot}_{\calH^{0,0}_2(\bz)}$ induces an identification 
	of the 	dual space $\calH^{0,\gamma_{p}}_{p}(\bz)'$ with 	
	$\calH^{0,\gamma_{p'}}_{p'}(\bz)$ and since 	
	$\ci(\bz)\subset\calH^{0,\gamma_{p'}}_{p'}(\bz)$, the result follows 
	if we 	can show that $G(\lambda)$ has an integral kernel 
	$k(\lambda)\in         
	\calS(\tilde{\Lambda}_{\theta},\calD(A_{p})\pit\ci(\bz))$. 
			A straightforward calculation shows that
			$G(\lambda)$ is a linear 	
combination of operators of the form 	
\begin{itemize}
	    \item[i)] $\varphi (\lambda - \widehat{A}_{p})^{-1} \psi$,
	              where $\varphi,\psi \in C^\infty([0,1[)$ have disjoint support 
		            and either $\varphi$ or $\psi$ is a cut-off function; 
	    \item[ii)] $\varphi_{0} (\lambda - \Delta_{2\bz}^p)^{-1} \psi_{0}$,
                 where $\varphi_{0},\psi_{0}\in C^\infty(]0,1[)$ have disjoint support; 
	    \item[iii)] $\varphi_{1} \{ (\lambda - \widehat{A}_{p})^{-1}
	                - (\lambda - \Delta_{2\bz}^p)^{-1} \} \psi_{1}$,
	                where $\varphi_{1},\psi_{1}\in C^\infty(]\varrho,1[)$. 
	\end{itemize}
	Both $(\lambda - \widehat{A}_{p})^{-1}$ and $(\lambda - 
	\Delta_{2\bz})^{-1}$ are, in particular, parameter-dependent 
	pseudodifferential parametrices on $]\varrho , 1 [ \times X$ to the 
	same operator $\lambda - \Delta$. Hence they coincide modulo smoothing 
	operators, and the terms from iii) are integral operators with a 
	parameter-dependent kernel belonging to 
	$\calS(\tilde{\Lambda}_{\theta}, C^{\infty,\infty}(\bz) 
	\pit C^{\infty,\infty}(\bz))$ where 
	$C^{\infty,\infty}(\bz)$ denotes the space of smooth functions on 
	$\bz$ vanishing of infinite order at the boundary. The same is true 
	for the terms from ii) due to the disjoint support of $\varphi_{0}$ 
	and $\psi_{0}$. Clearly, such integral operators have the required 
	property. 		
	It remains to consider terms from i). Since they are located 
	near the boundary, we can describe their kernel in the 
	splitting of coordinates $(t,\vartheta)$. It is given by
	$$k(\lambda,t,\vartheta,t^\prime,\vartheta^\prime)= \varphi(t) 
	\psi(t^\prime) \int_{\rz^{2}} 	      
	e^{i((t,\vartheta)-(t^\prime,\vartheta^\prime))\eta}
	\frac{1}{\lambda - |\eta|^{2}} \dbar \eta,$$	
	where, for abbreviation, we write $(t,\vartheta) := (t \cos \vartheta, 
	t\sin \vartheta)$.
	If $\psi$ is a cut-off function (hence $\varphi\in\cicomp(]0,1[)$), 
	this kernel belongs to 	
	$\calS(\tilde{\Lambda}_\theta,\calC^{\infty,\infty}(\bz)
	\pit\ci(\bz))$: indeed,  
	$|(t,\vartheta)-(t^\prime,\vartheta^\prime)| \ge c > 0$, in view of 
	the fact that the supports of $\varphi$ and $\psi$ are disjoint.
	Next suppose $\varphi$ is a cut-off function and $\psi \in 
	\cicomp({]0,1[})$. A Taylor expansion in $t$ of the 
	integral shows that $k=k_{0}+k_{1}+k_{2}$, where	    
	$$k_{0}(\lambda,t,\vartheta,t^\prime,\vartheta^\prime)=              
	\varphi(t) \psi(t^\prime) \int 	      
	e^{-i(t^\prime,\vartheta^\prime) \eta}
	\frac{1}{\lambda - |\eta|^{2}} \dbar \eta
	\in \calS(\tilde{\Lambda}_{\theta},
	\mathrm{span}(\varphi)\otimes\ci(\bz)),$$
	\begin{eqnarray*}
	k_{1}(\lambda,t,\vartheta,t^\prime,\vartheta^\prime)
	\!\! & \!\! = \!\! & \!\!
	\varphi(t) \psi(t^\prime) \, t
	\, \int e^{-i(t^\prime,\vartheta^\prime) \eta}
	(\eta_{1} \cos \vartheta + \eta_{2} \sin \vartheta)
	\frac{1}{\lambda - |\eta|^{2}} \dbar \eta
	\\
	\!\! & \!\! \in \!\! & \!\!
	\calS(\tilde{\Lambda}_{\theta},
	\mathrm{span}(\varphi\, t e^{i\vartheta},\varphi\, t e^{-i\vartheta})
	\otimes\ci(\bz)),
	\end{eqnarray*}
	while
	\[
	k_{2}(\lambda,t,\vartheta,t^\prime,\vartheta^\prime)
	=
	\varphi(t) \psi(t^\prime) \, t^{2}
	\, \int_{0}^{1} \int e^{i((st,\vartheta)-(t^\prime,\vartheta^\prime)) \eta}
	(\eta_{1} \cos \vartheta + \eta_{2} \sin \vartheta)^{2}
	\frac{1}{\lambda - |\eta|^{2}} (1-s) \dbar \eta ds.	
	\]
	The fact that the supports of $\varphi$ and $\psi$ are disjoint
	shows that $| (st,\vartheta)-(t^\prime,\vartheta^\prime) |$
	is bounded away from $0$ uniformly in $s$, hence
	$k_{2} \in \calS(\tilde{\Lambda}_{\theta}, \varphi t^{2} \ci(\bz)
	\pit \ci(\bz))$.
	
	Since $\mathrm{span}(\varphi, \varphi t e^{i\vartheta},
	\varphi t e^{-i\vartheta}) \subset \calD(\amin^p)$ and 
	$\varphi t^{2} \ci(\bz) \subset 
	\calH^{\infty,\gamma_{p}+2}_{p}(\bz) \subset \calD(\amin^p)$,
	we conclude that $k \in \calS(\tilde{\Lambda}_{\theta}, \calD(A_{p}) 
	\pit \ci(\bz))$ also in this case.
    \end{proof}
    
    \begin{proposition}
	If $R_{p}(\lambda)$ is as in \eqref{Res1}, then 
	\[
	(\lambda - A_{p}) R_{p}(\lambda) - 1 \in
		\calS(\tilde{\Lambda}_{\theta},
	      \calL(\calH^{0,\gamma_{p}}_{p}(\bz))),
	\qquad
	R_{p}(\lambda)(\lambda - A_{p}) - 1 \in 
	\calS(\tilde{\Lambda}_{\theta},
	      \calL(\calD(A_{p}))).
	\]
	Here, $\tilde{\Lambda}_\theta$ is the truncated sector as in Lemma 
	\text{\rm\ref{cutoff}}. In particular, $(\lambda - A_{p})^{-1}$ 
	exists for $\lambda \in \tilde{\Lambda}_{\theta}$ with $|\lambda|$ 
	sufficiently large, and we have 
        $$\|(\lambda-A_p)^{-1}\|_{{\calL}(\calH^{0,\gamma_p}_p(\bz))}=
O(|\lambda|^{-1}).$$
    \end{proposition}
    \begin{proof}
	Let us show the first statement. To this end, write 
	$$\lambda-A_p=\sigma\,(\lambda+\widehat{\Delta}^p)\,\sigma_0+
	(1-\sigma)\,(\lambda+\Delta_{2\bz}^p)\,(1-\sigma_1)$$
	with cut-off functions $\sigma,\sigma_j\in\cicomp([0,1[)$ satisfying 
	$\sigma\sigma_1=\sigma_1$, $\sigma_0\sigma=\sigma$, and 
	$\sigma\equiv\sigma_j\equiv1$ on $[0,\varrho]$. Then 
	$$(\lambda-A_p)R_p(\lambda)=
	\sigma\,(\lambda+\widehat{\Delta}^p)\,\sigma_0\,R_p(\lambda)+
	(1-\sigma)\,(\lambda+\Delta_{2\bz}^p)\,(1-\sigma_1)\,R_p(\lambda).$$
	To treat the first summand on the right-hand side, choose a 
	representation of $R_p(\lambda)$ with cut-off function $\omega$ such 
	that $\omega\sigma_0=\sigma_0$. According to the previous Lemma 
	\ref{cutoff} and the fact that the operator norm of 
	$\sigma\,(\lambda+\widehat{\Delta}^p)\,\sigma_0:\calD(A_p)
	\to\calH^{0,\gamma_p}_p(\bz)$ 
	is $O(|\lambda|)$,  
	$$\sigma\,(\lambda+\widehat{\Delta}^p)\,\sigma_0\,R_p(\lambda)
	\equiv
	\sigma\,(\lambda+\widehat{\Delta}^p)\,\sigma_0
	\omega\,(\lambda+\widehat{\Delta}^p)^{-1}\,\omega_0
	=
	\sigma\,(\lambda+
	\widehat{\Delta}^p)(\lambda+\widehat{\Delta}^p)^{-1}\,\omega_0
	=\sigma$$
	modulo a remainder in 
	$\calS(\tilde{\Lambda}_\theta,\calL(\calH^{0,\gamma_p}_p(\bz)))$; 
	note that the factor $\sigma_0\omega$ can be omitted due to the 
	locality of $\widehat{\Delta}$. For the second summand choose 
	$\omega_0$ such that $\sigma_1\omega_0=\omega_0$ to obtain 
	analogously 
	$(1-\sigma)\,(\lambda+\Delta_{2\bz}^p)\,(1-\sigma_1)\,R_p(\lambda)\equiv 
	1-\sigma$ modulo a remainder of the same type. The proof of the second 
	statement is analogous.  
        Finally the norm estimate is immediate from the form of $R_p$ in Lemma
        \ref{cutoff}.
    \end{proof}
    
\begin{corollary}\label{holomsg}Let $\Delta$ be the
Laplace-Beltrami operator of Example {\rm \ref{laplace4}}.
Then $-\Delta$ is the generator of
a holomorphic semigroup. 
\end{corollary}

Note that this is already sufficient for 
the solution of certain semilinear evolution equations, cf.\ e.g. Pazy
\cite[Theorem 6.3.1]{Pazy}.

\section{Quasilinear parabolic equations}\label{semi5}
       
    In the previous section we saw that the boundedness of the purely imaginary powers 
    implies the solvability of associated parabolic initial 
    value problems and the maximal regularity of the 
    solution. In turn, the knowledge of maximal regularity is 
    important for the investigation of non-linear equations, as we 
    want to illustrate in this section. Following the concept of 
    Cl\'ement and Li \cite{ClLi}, we will consider examples of quasilinear
    evolution equations.
    
    Let $E=(E_{0}, E_{1})$ be a couple of Banach spaces such that 
    $E_{1}$ is 
    densely and continuously embedded into $E_{0}$. For $1 < q < \infty$ 
    denote by 
    \[
    E_{1-\frac{1}{q},q} := (E_{1},E_{0})_{\frac{1}{q},q}=
    (E_{0},E_{1})_{1-\frac{1}{q},q}
    \]
    the space given by the real interpolation method 
    $(.,.)_{\vartheta,q}$.

    Let $-P \in \calL(E_{1},E_{0})$ 
    be the infinitesimal generator of an analytic semigroup in 
    $E_{0}$ with $\calD (P) = E_{1}$. For $T>0$ and $f \in 
    L_q([0,T];E_{0})$, a function $u \in 
    W^1_q([0,T];E_{0}) \cap L_q([0,T];E_{1})$ is called a strict 
    solution of the problem
    \begin{equation}
	\label{para}
	\left\{
	\begin{array}{l}
	    \dot{u}(\tau) + Pu(\tau) = f(\tau) \mbox{ on $]0,T[$}
	    \\
	    u(0) = u_{0},
	\end{array}
	\right.
    \end{equation}
    if $u$ satisfies (\ref{para}) in the $L_q([0,T];E_{0})$ sense. 
    It is known that (\ref{para}) with $f \equiv 0$ has a strict 
    solution 
    if and only if $u_{0} \in E_{1-\frac{1}{q},q}$ (see, e.g., 
    \cite[Theorem 4.10.2]{Aman}).
    
    We will say that
    $P \in \calL(E_{1},E_{0})$ belongs to the class $\Mr(q,(E_0,E_1))$ if
    for every $f \in L_q([0,T];E_{0})$ and $u_{0} \in 
    E_{1-\frac{1}{q},q}$, 
    there exists a unique strict solution $u \in W^1_q([0,T];E_{0}) 
    \cap L_q([0,T];E_{1})$ of (\ref{para}) 
    and if there exists $M > 0$, independent of $f$ and $u_{0}$, such 
    that
    \[
      \int_{0}^{T} \| \dot{u}(\tau) \|^q_{E_{0}}d\tau
    + \int_{0}^{T} \| Pu(\tau) \|^q_{E_{0}} d\tau \le
      M \Big( \int_{0}^{T} \| f(\tau) \|_{E_{0}}^q d\tau 
               + 
	       \| u_{0} \|_{ E_{1-\frac{1}{q},q}}^q \Big).
    \]
    Cl\'ement and Li considered  the 
    quasilinear problem
 	\begin{equation}
	    \label{slpara}
	    \left\{
	    \begin{array}{l}
		\dot{u}(\tau) + A(u)u(\tau) = f(\tau,u(\tau)) + g(\tau) \mbox{ on $]0,T_{0}[$}
		\\
		u(0) = u_{0},
	    \end{array}
	    \right.
	\end{equation}
 where $T_0> 0$ and $A$, $f$, and $g$ are supposed to satisfy the following assumptions:
\begin{itemize}
\item[(H1)]$A\in {\mathcal C}^{1-}(U,{\mathcal L}(E_1, E_0))$ 
for some open neighborhood $U$ of $u_0$ in $E_{1-\frac1q,q}$, and $A(u_0)\in
\Mr(q,(E_0,E_1))$;
\item [(H2)]$f\in {\mathcal C}^{1-,1-}([0,T_0]\times U,E_0)$;
\item [(H3)]$g\in L_q([0,T_0], E_0)$.
\end{itemize}

Their main result then is:  
    \begin{theorem}
	\label{ClLi}
	Under hypotheses  {\rm (H1)}, {\rm (H2)}, and {\rm (H3)} there exists a
	$T_1\in ]0,T_0]$ 
  and a unique function $u\in W^1_q([0,T_{1}];E_{0}) \cap L_q([0,T_{1}];E_{1})
  \cap {\mathcal C}([0,T_1],E_{1-\frac1q,q})$
  satisfying \eqref{slpara} on $]0,T_1[$.
\end{theorem}
    
We shall now show how this theorem can be applied to certain equations on 
manifolds with conical singularities.
To this end we shall verify the conditions for some operators $A$ 
and functions $f$. In the following, we will fix 
$E_0=\calH^{0,\gamma_{p}}_{q}(\bz)$ and $
E_1=\calH^{2,\gamma_{p}+2}_{q}(\bz)$ with $1 < p,q < \infty$.

\subsection{A Lipschitz continuous family of Laplace type operators}

    By Example \ref{exlaplmaxreg}, the operator $-\Delta + c$ on 
    $\bz$ 
    admits imaginary powers such that $\| (-\Delta +c)^{iy}
    \|_{\calL(\calH^{0,\gamma_{p}}_{q}(\bz))}$ $\le C e^{\theta |y|}$
    with $0 < \theta < \frac{\pi}{2}$ if $c$ is large enough
    and $2\max(p,p^\prime)-1<n$. 
    We apply a slight extension of 
    the Dore-Venni Theorem \ref{dorevenni} valid for 
    arbitrary initial data $u_{0}$, see \cite{Aman}, Theorem III.4.10.7
    and conclude that
    $-\Delta \in \Mr(q,(\calH^{0,\gamma_{p}}_{q}(\bz),
    \calH^{2,\gamma_{p}+2}_{q}(\bz)))$
    for any $1 < p,q < \infty$.

The situation does not change very much if we replace $-\Delta$ by
$-b\Delta$, 
where $b$ is a smooth positive function on $\bz$ which is constant
at the boundary $\{t=0\}$, say $b|_{\{t=0\}} = b_0$: 
the principal symbol and the rescaled symbol of 
$-b\Delta$ are $-b\sigma_\psi^2(\Delta)$ and $-b_0\tilde\sigma_\psi^2(\Delta)$,
respectively.
They are invertible in the same sector as $-\sigma_\psi^2(\Delta)$  and 
 $-\tilde\sigma_\psi^2(\Delta)$, respectively. As condition (E)  holds for $-\Delta$, 
it also holds for $-b\Delta$.

Also the two model cone operators differ only by the constant $b_0$:
 $$-\widehat{b\Delta}= -b_0\widehat{\Delta}.$$ 
Hence $(\mathrm E_{\mathrm {min}})$ holds for $-\widehat{b\Delta}$ in the same sector it holds for
$-\widehat{\Delta}$. We may therefore apply Theorem \ref{criteria} and obtain:

\begin{proposition}\label{bDelta}
Given a smooth function $b$ on $\bz$ which is constant at $\partial \bz$, 
the operator $-b\Delta$ is an element of 
$\Mr(q,(\calH^{0,\gamma_{p}}_{q}(\bz),\calH^{2,\gamma_{p}+2}_{q}(\bz)))$
    for any $1 < p,q < \infty$ with 
    $2 \max(p,p^\prime)-1<n$.
\end{proposition}

We note the following simple lemma:

\begin{lemma}\label{embedding}For
$s>\frac{n+1}q$ and $\gamma\ge \frac{n+1}2$ we have 
$\calH^{s,\gamma}_q(\bz)\hookrightarrow 
{\mathcal C}_b(\bz)$, the space of bounded continuous functions on $\bz$.
\end{lemma}

\begin{proof}Outside a neighborhood of the boundary, the space
$\calH^{s,\gamma}_q (\bz)$ coincides with the standard Sobolev space
$H^{s}_q(\bz)$, so that our statement follows from the well-known embedding 
theorem. For functions in $\calH^{s,\gamma}_q (\bz)$  supported near 
$\partial \bz$, we apply the mapping $S_\gamma$ defined in \eqref{sgamma}.                                
\end{proof}

   Our next step is to study $E_{1-\frac1q,q}$. A precise description of
   this interpolation space requires the introduction of weighted
   Besov spaces on $\bz$. For our purposes, however, the following embedding
   statement is sufficient.
    
    \begin{lemma}
	\label{interp}
	Let $s_{0}$, $s_{1}$, $\gamma_{0}$, $\gamma_{1} \in \rz$,
	$0 < \vartheta < 1$ and $1 < q < \infty$. Then, for arbitrary
	$\delta, \varepsilon > 0$,
	\[
	( \calH_{q}^{s_{0},\gamma_{0}}(\bz), 
	\calH^{s_{1},\gamma_{1}}_{q}(\bz)
	)_{\vartheta,q} \hookrightarrow 
	\left\{
	\begin{array}{ll}
	    \calH_{q}^{s,\gamma-\varepsilon}(\bz)        & 
	    \mbox{if $q \le 2$}
	    \\
	    \calH_{q}^{s-\delta,\gamma-\varepsilon}(\bz) &
	    \mbox{if $q > 2$}
	\end{array}
	\right.
	\]
	with $s=(1-\vartheta)s_{0}+\vartheta s_{1}$, $\gamma = 
	(1-\vartheta)\gamma_{0}+\vartheta \gamma_{1}$.
    \end{lemma}
    \begin{proof}
	By definition of the cone Sobolev spaces, cf. \eqref{norm},
	the statement is true if we can show the following interpolation 
	result 
	for the local spaces (where for notational simplicity we suppress 
	writing $\rz^{1+n}$):
	\begin{equation}
	    \label{interA}
	    (H^{s_0,\gamma_0}_{q},H^{s_1,\gamma_1}_{q})_{\vartheta,q}
	    \hookrightarrow
	    \begin{cases}
		H^{s,\gamma-\varepsilon}_q & \mbox{ if }q\le 2 \\
		H^{s-\delta,\gamma-\varepsilon}_q & \mbox{ if }q>2,
	    \end{cases}
	\end{equation}
	where the $H^{r,\sigma}_q$ denotes the weighted space
	$e^{-\sigma\spk{t}}H^r_q(\rz^{1+n}_{(t,x)})$. Note that for
	$\gamma_0=\gamma_1=0$  
	\begin{equation}\label{interB}
	    (H^{s_0}_{q},H^{s_1}_{q})_{\vartheta,q}=B^s_{qq}
	\end{equation}
	is a Besov space; in this case \eqref{interA} follows from standard
	embedding properties (even $\varepsilon=0$ is true), cf. Triebel
	\cite{Trieb}. 

	To prove the general case we need to introduce some notation.
	The following method appears, e.g., in \cite{Trieb}, \cite{BeLo}.
	For a Banach space $Y$ and real $r$ we let $l^r_q(Y)$ denote 
	the space of all sequences $(y_k)_{k\in\gz}$ in $Y$ such that 
	$$\|(y_k)\|:=\Big(\smsum_{k\in\gz}(e^{r|k|}\|y_k\|_Y)^q
	\Big)^{1/q}<\infty.$$
	Then, if interpolation of $Y_0$ and $Y_1$ makes sense, i.e., 
	$(Y_{0},Y_{1})$ is an interpolation couple, by 
	Theorem 5.6.2 of \cite{BeLo}
	\begin{equation}\label{interC}
	    (l^{\gamma_0}_q(Y_{0}),l^{\gamma_1}_q(Y_{1}))_{\vartheta,q}=
	    l^\gamma_q((Y_0,Y_1)_{\vartheta,q}).
	\end{equation}
	Furthermore let us fix a function 
	$\varphi=\varphi(t)\in\cicomp(\rz)$ 
	supported in $[-1,1]$ and strictly positive in $]-1,1[$ such that
	$\sum\limits_{k\in\gz}\varphi(\cdot-k)\equiv1$ on all of $\rz$. 
	Then define
	$\varphi_k\in\ci(\rz^{1+n})$ by $\varphi_k(t,x)=\varphi(t-k)$. 
  For $u\in H^{r,\sigma}_q$ we can estimate 
	$$\|\varphi_k u\|_{H^r_q}=
     \|(\varphi_k e^{-\sigma\spk{t}})(e^{\sigma\spk{t}}u)\|_{H^r_q}\le
	   C\sup_{k-1\le t\le k+1}e^{-\sigma\spk{t}}\|u\|_{H^{r,\sigma}_q}\le
	   Ce^{-\sigma|k|}\|u\|_{H^{r,\sigma}_q},$$
  	with a constant $C$ independent of $k$.
	Then we use the fact that the operator norm of a map $u\mapsto au$ in
	$H^r_q$ for some $a\in\calC^\infty_b(\rz^{1+n})$ can
	be estimated by finitely many terms $\|D^\alpha a\|_\infty$.
	Hence, for any $\sigma'<\sigma$, the map 
	\begin{equation}\label{interD}
	    S:u\mapsto(\varphi_k u)_{k\in\gz},\quad H^{r,\sigma}_q\to
	    l^{\sigma'}_q(H^r_q)
	\end{equation}
	is well-defined and continuous. On the other hand, if 
	$(u_k)_{k\in\gz}\in l^{\sigma'}_q(H^r_q)$ is given and
	$\psi\in\cicomp(\rz)$ is chosen in such a way that 
	$\psi\varphi=\varphi$,
	and we set $\psi_k(t,x)=\psi(t-k)$ then 
	$$\|\psi_k u_{k} \|_{H^{r,\sigma''}_q}=\|(\psi_k 
	e^{\sigma''\spk{t}})u_k\|_{H^r_q}\le
	Ce^{\sigma''|k|}\|u_k\|_{H^r_q}$$
	by an argument analogous to the above one. This together with 
	H\"older's 
	inequality shows that for any $\sigma''<\sigma'$ the map 
	\begin{equation}\label{interE}
	    R:(u_k)_{k\in\gz}\mapsto\smsum_{k\in\gz}\psi_k u_k,\quad 
	    l^{\sigma'}_q(H^r_q)\to     H^{r,\sigma''}_q 
	\end{equation}
	is well-defined and continuous. Clearly $RSu=u$ 
	for any $u\in H^{r,\sigma}_q$, by the choice of $\psi$. 
	From \eqref{interD}, \eqref{interC}, and \eqref{interB} we now 
	obtain 
	$$S:(H^{s_0,\gamma_0}_{q},H^{s_1,\gamma_1}_{q})_{\vartheta,q}\to 
	l^{\gamma-\varepsilon/2}(B^s_{qq})\hookrightarrow
	l^{\gamma-\varepsilon/2}(H^{s-\delta}_q),$$
	where $\delta=0$ if $q\le2$, $\delta>0$ if $q>2$, and 
	$\varepsilon>0$.
	Applying \eqref{interE}
	we get that 
	$$\iota=RS:(H^{s_0,\gamma_0}_{q},H^{s_1,\gamma_1}_{q})_{\vartheta,q}\to 
	H^{s-\delta,\gamma-\varepsilon}_q$$ 
        with $\delta,\varepsilon > 0$ arbitrarily small is a 
	continuous embedding. 
    \end{proof}
    
    \begin{corollary}
	\label{interp2}
	$( \calH_{q}^{2,\gamma_{p}+2}(\bz), \calH^{0,\gamma_{p}}_{q}(\bz) 
	)_{\frac{1}{q},q}
	\hookrightarrow \calH^{s,\delta}_{q}(\bz)$ for 
	any $s < \sfrac{2}{q^\prime}$ and any $\delta < \gamma_{p} + 
	\sfrac{2}{q^\prime}$.
    \end{corollary}
In the sequel, we shall denote by $t$ a smooth, strictly positive function on 
$\bz$ which coincides with the distance to the boundary 
(i.e. the coordinate $t$ employed above) in a collar neighborhood.

\begin{lemma}\label{embedinC}
    Let $c >0$ and $1<p,q<\infty$ with 
$p\ge \frac{n+1}{2+c}$,
$q>\max\left(\frac{n+3}2,\frac{2p}{(2+c) p-(n+1)}\right)$.
Then  $$E_{1-\frac1q,q}\hookrightarrow t^{-c}{\mathcal C}_b(\bz).$$
\end{lemma}

\begin{proof}The  conditions on $p$ and $q$ allow us to find $s,\delta$ with
$ \frac{(n+1)}q<s<\frac{2}{q'}$ and 
$\frac{n+1}{2}-c\le\delta<\gamma_p+\frac{2}{q'}.$
According to Corollary \ref{interp2} and Lemma \ref{embedding}
we have 
$E_{1-\frac1q,q}
\hookrightarrow \calH^{s,\delta}_{q}(\bz)
\hookrightarrow t^{-c}\calH^{s,\frac{n+1}2}_{q}(\bz)\hookrightarrow
t^{-c}
{\mathcal C}_b(\bz).$
Note that the second inclusion is immediate from Definition \ref{sobolev} with interpolation.
\end{proof}

Not for every $n$ it is possible to find $p$ and $q$ satisfying the 
hypotheses of Lemma \ref{embedinC} and the inequality $2 
\max(p,p^\prime) - 1 < n$. 
However, all these requirements can be 
fulfilled at the same time when
$n\ge 4$, i.e., dim\,$\bz\ge 5$.
\begin{theorem}\label{aDelta} 
Let $n\ge 4$ and $c>0$.
Choose $1<p,q<\infty$ with $2\max(p,p')<{n+1}$, 
$p\ge \frac{n+1}{2+c}$, and
$q>\max\left(\frac{n+3}2,\frac{2p}{(2+c) p-(n+1)}\right)$. 
Fix a smooth initial value  $u_0$ which vanishes to infinite order at $\partial \bz$
and a bounded neighborhood $U$ of $u_0$ in $E_{1-\frac1q,q}$.
Let $a$ be a smooth, strictly positive function on $\cz\cong \rz^2$. Then
the operator function $A(u) = -a(t^c u)\Delta$, $u\in U$, satisfies {\rm(H1)}. 
\end{theorem}

\begin{proof}
The function $t^c u_0$ is smooth on $\bz$, hence so is $b=a(t^c u_0)$.
In addition, $b$ is positive and constant at the boundary.
By Proposition
\ref{bDelta}, $A(u_0) = -b\Delta$ belongs to  
$\Mr(q,(\calH^{0,\gamma_{p}}_{q}(\bz),\calH^{2,\gamma_{p}+2}_{q}(\bz)))$. 

As $u$ varies over a bounded neighborhood $U$ of $u_0$ in $E_{1-\frac1q,q}$, 
the functions
$t^c u$ vary over a bounded set in ${\mathcal C}_b(\bz)$.  In
particular, $a(t^c u)$ is a continuous, bounded, and strictly
positive function on $\bz$. Hence $A(u)$ is an element of $\calL(E_1,E_0)$
for each $u$. Moreover, 
\begin{eqnarray*}
\|A(u_1)-A(u_2)\|_{\calL(E_1,E_0)}&\le&
\|a(t^c u_1)
-a(t^c u_2)\|_{L_\infty(\bz)}\|\Delta\|_{\calL(E_1,E_0)}
\\
&\le&
C
\|t^c u_1-t^c u_2\|_{\calC_b(\bz)}\|\Delta\|_{\calL(E_1,E_0)}
\le C \|u_1-u_2\|_{E_{1-\frac1q,q}} \|\Delta\|_{\calL(E_1,E_0)},
\end{eqnarray*}
where $C$ is the maximum of $|a'(s)|$ as $s$ varies over the bounded set of all
values of $t^c u$, $u\in U$. 
\end{proof} 
\begin{remark}\label{integer}
In case $c \in \nz$, the initial
value $u_0$ can be chosen in ${\mathcal C}^\infty(\bz)$.
Theorem \ref{aDelta} 
extends to the case where $c$
is a smooth real-valued
function on $\bz$ which is positive and constant at the boundary.

%
\end{remark}


\subsection{Lipschitz continuity of the functions 
\boldmath{$|u|^\alpha$}}

  Let us now try to find functions $f$ satisfying hypothesis (H2).
 Here is a first simple example.
    
  \begin{example}
	\label{experiodic}
	{\rm (H2)} holds for $f(\tau,u) = h(u)$, with
	$h \in \calC(\rz)$ such that $h(0)=0$ and $|h(s)-h(s^\prime)| \le 
	M|s-s^\prime|$ for some $M \ge 0$, uniformly in $s,s^\prime \in \rz$.
    \end{example}
    %
    This follows from the observation that the mapping
    \[
    u \mapsto h(u(.)) \, : \, L_q(\Omega) \rightarrow L_q(\Omega)
    \]
    is Lipschitz continuous for any measure space $\Omega$, in 
    particular for $\calH^{0,\gamma_{p}}_{q}(\bz)$, which is a 
    weighted $L_q$-space on $\bz$.
    
    As mentioned in the introduction, 
    nonlinearities of the type $|u|^\alpha$ or $u^\alpha$ are relevant for 
    applications.
    It is then interesting to find out whether  
    a term of this kind fulfills (H2). 
    We shall show the following:
       
    \begin{theorem}
	\label{nonlinear}
The function $f(\tau,u)=|u|^\alpha$ satisfies {\rm (H2)}
	for all $1\le\alpha<\alpha^{*}$, where $\alpha^{*}$ is determined as follows: 
  \begin{itemize}
   \item[a)] If $\frac{2p}{q^\prime}<n+1$ then \ 
     $\alpha^{*}=\begin{cases}
                   \frac{n+1}{n+1-2p/q^\prime} & q\ge\frac{n+3}{2}\\
                   \min\left(\frac{n+1}{n+1-2p/q^\prime},
                             \frac{n+1}{n+1-2q/q^\prime}\right) & 
                       q<\frac{n+3}{2}
                 \end{cases},$
   \item[b)] If $\frac{2p}{q^\prime}\ge n+1$ then \ 
     $\alpha^{*}=\begin{cases}
                   \infty & q\ge\frac{n+3}{2}\\
                   \frac{n+1}{n+1-2q/q^\prime} & q<\frac{n+3}{2}
                 \end{cases}.$
  \end{itemize}
 \end{theorem}

    \begin{corollary}
	\label{cornonlinear}
	Let $n\ge 4$ and $f(t,u) = |u|^\alpha$. 
\begin{itemize} 
	
\item [{\rm a)}] Hypothesis  {\rm (H2)} is satisfied for arbitrary $\alpha \ge 1$, 
provided we choose $p<\frac{n+1}2$ sufficiently close to $\frac{n+1}2$ and $q\in ]1,\infty[$ sufficiently large.
In this case, $2\max(p,p')-1<n$, and $p$ and $q$ satisfy also
the conditions of Theorem {\rm \ref{aDelta}}.

\item [{\rm b)}] Given $p$ with $2\max(p,p')-1<n$,  {\rm (H2)} 
holds for  $1 \le \alpha < \sfrac{n+1}{n+3-2p}$ with $q=p$.

\item [{\rm c)}] Hypothesis  {\rm (H2)} is satisfied for $1\le\alpha<\frac{n+1}2$
if  $q=p<\frac{n+1}2$ is sufficiently close to $\frac{n+1}2$.
\end{itemize}
    \end{corollary}
  \begin{proof} 
For $p<\frac{n+1}2$ sufficiently close to $\frac{n+1}2$, we have $2\max(p,p^\prime)-1<n$. 
Conversely, the condition
$2\max(p,p^\prime)-1<n$ implies $2p<n+1$. So the assertions follow
from  Theorem \ref{nonlinear}a). 
  \end{proof}

    \begin{remark}\label{u^alpha}
	The statements of Theorem {\rm\ref{nonlinear}} and Corollary 
	{\rm\ref{cornonlinear}} remain true in case $f(t,u) = u^\alpha$ 
	with a natural number $\alpha$  satisfying the corresponding 
	conditions.
    \end{remark}
    We shall prove Theorem \ref{nonlinear} after Theorem \ref{lipF2}. 
    We first note that, as a consequence of
    Corollary \ref{interp2}, the Lipschitz continuity of a map on bounded subsets of 
    $E_{1-\frac{1}{q},q}$ follows from its Lipschitz continuity 
    on bounded subsets of $\calH^{s,\delta}_{q}(\bz)$. Let us now 
    introduce $\calH^{s,\gamma}_{q}(X^\wedge)$ as the space of all 
    distributions $u$ on $X^\wedge = \rp \times X$ such that 
    $S_{\gamma}u \in H^s_{q}(\rz \times X)$, where 
    $S_{\gamma}$ is the map introduced in \eqref{sgamma} (with $X$ 
    instead of $\rz^n$); the norm is given by 
    $\| u \|_{\calH^{s,\gamma}_{q}(X^\wedge)} = 
     \| S_{\gamma} u \|_{H^s_{q}(\rz \times X)}$.
     
    Moreover, we denote by   
    $\calH^{s,\gamma}_{q}(X^\wedge)_{0}$ the subspace of all 
    $u \in \calH^{s,\gamma}_{q}(X^\wedge)$
    supported in $[0,1[ \times X$.

    \begin{lemma}
	\label{badembed}
	Let $1 < q \le \tilde{q} < \infty$ and $\eps > 0$. Then
	$\calH_{\tilde{q}}^{0,\gamma+\eps}(X^\wedge)_{0}
	\hookrightarrow\calH^{0,\gamma}_{q}(X^\wedge)_{0}$.
    \end{lemma}
    \begin{proof}
	By the definition of $\calH^{0,\gamma+\eps}_{\tilde{q}}(X^\wedge)$ 
	we have
	\begin{align*}
	    \| u \|^q_{\calH^{0,\gamma}_{q}(X^\wedge)}
	    &  =  
	           \int_{\rp \times X} 
	                r^{\left( \frac{n+1}{2} - \gamma \right) q} \,
			|u(r,x)|^q \, \frac{dr}{r} dx	    
	    \\
	    &  =  
	           \int_{]0,1] \times X} 
	                r^{\frac{q}{\tilde{q}} - 1 + \eps q} \,
	           \big(r^{\left( \frac{n+1}{2} - \gamma - 
		   \frac{1}{\tilde{q}}
			           - \eps \right) q} \, |u(r,x)|^q\big)\, dr dx. 
	\end{align*}
	The second factor of the integrand belongs to 
	$L_{\frac{\tilde{q}}{q}}(X^\wedge,drdx)$. 
	Using H\"older's inequality we can estimate the integral from above 
	by  
	$$\Big(\int_{]0,1] \times X} 
	       r^{ - 1 + \eps \frac{q\tilde{q}}{\tilde{q}-q}}drdx 
	       \Big)^\frac{\tilde{q}-q}{\tilde{q}}
	   \Big(\int_{]0,1] \times X}
	        r^{\left( \frac{n+1}{2} - \gamma - \eps \right) \tilde{q}}\,
	        |u(r,x)|^{\tilde{q}}\,\frac{dr}{r}dx 
		\Big)^{\frac{q}{\tilde{q}}}
	   \le C\,\|u\|^q_{\calH^{0,\gamma+\eps}_{\tilde{q}}(X^\wedge)}.$$
	   This shows the continuity of the stated embedding. 
    \end{proof}
    

        \begin{lemma}
	\label{contF}
	Let $\alpha \ge 1$,  $s \ge 0$, and $1 < q < \infty$.

        For $sq < n+1$,  $u \mapsto |u|^\alpha$ maps bounded subsets of
	$\calH^{s,\varrho}_{q}(X^\wedge)$ to bounded subsets of
	$\calH^{0,\tilde{\varrho}}_{\tilde{q}}(X^\wedge)$
        for $\tilde{\varrho}=\alpha\varrho-(\alpha-1)
	\frac{n+1}{2}$ and any $\tilde q$ in $[\frac q\alpha, \frac
        q\alpha\frac{n+1}{n+1-sq}]\ \cap\  ]1,\infty[$.

        For $sq \ge n+1$, the same is true  for $\tilde \varrho$ as before
        and $\tilde q>\max(1,q/\alpha)$.

	In both cases there is a positive constant $C$ such that 
	\[
	\| \, |u|^\alpha \, 
	\|_{\calH^{0,\tilde{\varrho}}_{\tilde{q}}(X^\wedge)}
	\le
	C \| u \|^\alpha_{\calH^{s,\varrho}_{q}(X^\wedge)}
	\hspace{1cm} \forall u \in \calH^{s,\varrho}_{q}(X^\wedge).
	\]
    \end{lemma}

    \begin{proof}
        It is well-known that $H^s_p(\rz^n)\hookrightarrow H^t_q(\rz^n)$
provided
        $1<p\le q<\infty$ and $s-n/p\ge t-n/q$,
cf. \cite[2.8.1 Remark 2]{Trieb}.
        Hence we have  $H^{s}_{q}(\rz^{1+n})
\hookrightarrow
	L_{\alpha \tilde{q}}(\rz^{1+n})$ and
	\begin{align*}
	    \| \, |u|^\alpha \, 
	    \|_{\calH^{0,\tilde{\varrho}}_{\tilde{q}}(X^\wedge)}
	    & = 
	    \| \, S_{\tilde{\varrho}} |u|^\alpha \, 
	    \|_{L_{\tilde{q}}(\rz \times X)}
	    =
	    \| \, |S_{\varrho} u|^\alpha \,
	    \|_{L_{\tilde{q}}(\rz \times X)}\\
	    & =
	    \| S_{\varrho} u 
	    \|^\alpha_{L_{\alpha \tilde{q}}(\rz \times X)}
	    \le 
	    C \| S_{\varrho} u \|^\alpha_{H^s_{q}(\rz \times X)}
	    = 
	    C \| u \|^\alpha_{\calH^{s,\varrho}_{q}(X^\wedge)}.	    
	\end{align*}
    \end{proof}
    
 
\begin{theorem}
	\label{lipF2}
	Let $1 < q < \infty$ and $\gamma<\delta$. If $0<s<\frac{n+1}{q}$, 
	then the map 
	$$u\mapsto|u|^\alpha:\calH^{s,\delta}_q(\bz)\to
	\calH^{0,\gamma}_q(\bz)$$
	is Lipschitz continuous on bounded sets, whenever 
	\begin{align}
	    1 & \le \alpha  <  
           \min\Big(
		            \frac{n+1-2\gamma}{n+1-2\delta},\frac{n+1}{n+1-sq} 
			        \Big) 
        & \text{if }\delta<\frac{n+1}{2},\label{condalpha1}\\
	    1 & \le \alpha < \frac{n+1}{n+1-sq} 
			 & \text{if }\delta\ge\frac{n+1}{2}\label{condalpha2}.
	\end{align}
For $sq\ge n+1$, the same result is true; the upper bound for $\alpha $ then is
$\frac{n+1-2\gamma}{n+1-2\delta}$ in case $\delta <\frac{n+1}2$
and $\infty$ for
$\delta \ge \frac{n+1}2$.

    \end{theorem}

  \begin{proof}
We start with the simple observation that, for every measure space $\Omega$ and 
for every choice of $r$ such that $1<r\le\alpha r<\infty$, the map
	\[
	v \mapsto |v|^\alpha \; : \; L_{\alpha r}(\Omega)
	\rightarrow L_{r}(\Omega)
	\]
	is Lipschitz continuous on bounded sets. Indeed, this is a 
	straightforward consequence of H\"older's inequality and the fact 
	that $|x^\alpha - y^\alpha| \le \alpha \max\{x^{\alpha-1}, 
	y^{\alpha-1}\} |x-y|$ for any $x,y \ge 0$ and $\alpha \ge 1$.

	The crucial part of the proof concerns the analysis near the boundary,
        i.e., the Lipschitz continuity of the map
	\begin{equation}\label{lip}
   u\mapsto|u|^\alpha:\calH^{s,\delta}_q(X^\wedge)_0\longrightarrow
	    \calH^{0,\gamma}_q(X^\wedge)_0.
	\end{equation} 
	Assume we have proved this. 
        Combining the above observation with the embedding 
        $H^s_q(2\bz)\hookrightarrow L_{\alpha q}(2\bz)$ which is valid for all
        $\alpha \ge 1$ in case $sq\ge n+1$ and for
	$1 \le \alpha < \sfrac{n+1}{n+1-sq}$ in case $sq<n+1$,
        we immediately see that the map $u \mapsto |u|^\alpha$ is Lipschitz
        continuous from 
	bounded subsets of $H^s_{q}(2\bz)$ to $L_q(2\bz)$.
	Now, for an arbitrary cut-off function $\omega$,  
	choose cut-off functions $\sigma_{1}$ and $\sigma_{2}$ such that
	$\sigma_{1} \equiv 1$ on the support of 
	$\omega$ and $\omega \equiv 1$ on the support of $\sigma_{2}$.
	Then, by the definition of $\calH^{0,\gamma}_{q}(\bz)$, we have, 
	for $u$ and $v$ running through bounded subsets,
	\begin{align*}
		  \| \; |u|^\alpha - |v|^\alpha \; 
		  \|_{\calH^{0,\gamma}_{q}(\bz)}
		  & \le
		  \| \omega(|u|^\alpha - 
		  |v|^\alpha) \|_{\calH^{0,\gamma}_{q}(X^\wedge)}
		  + \|(1- \omega)(|u|^\alpha - 
		  |v|^\alpha) \|_{L_{q}(2\bz)}
		  \\
		  & \le 
		  \| \sigma_{1}^\alpha(|u|^\alpha - 
		  |v|^\alpha) \|_{\calH^{0,\gamma}_{q}(X^\wedge)}
		  + \|(1- \sigma_{2})^\alpha(|v|^\alpha - 
		  |v|^\alpha) \|_{L_{q}(2\bz)}
		  \\
		  & \le 
		  C \big( \| \sigma_{1}(u - v) \|_{\calH^{s,\delta}_{q}(X^\wedge)}
				 +
				 \| (1-\sigma_{2}) (u - v) \|_{H^s_{q}(2\bz)} \big)
       \le C \| u - v \|_{\calH^{s,\delta}_{q}(\bz)}.
	\end{align*}
	We next verify \eqref{lip}. In case $sq<n+1$, we set 
	$\tilde{q} = \sfrac{q}{\alpha} \, \sfrac{n+1}{n+1-sq}$; for 
        $sq\ge n+1$ we choose $\tilde q>q $ arbitrary. Note that also in the
        first case, our assumption on $\alpha$ implies $\tilde q>q$. 
        For arbitrary $\beta > 0$, we obtain 
	\begin{equation}
	    \label{eq:lip1}
	    \begin{array}{l}
		\| \; |u|^\alpha - |v|^\alpha \; 
		\|_{\calH^{0,\gamma}_{q}(X^\wedge)}
		\le
		C \| \; |u|^\alpha - |v|^\alpha \;
		  \|_{\calH^{0,\gamma+\beta}_{\tilde{q}}(X^\wedge)} =
		\\
		\hspace{1cm}
		= C \| S_{\gamma+\beta}(|u|^\alpha) - 
		       S_{\gamma+\beta}(|v|^\alpha) 
		    \|_{L_{\tilde{q}}(\rz\times X)} 
		= C \| \;
		       |S_{\tilde{\gamma}}u|^\alpha - 
		       |S_{\tilde{\gamma}}v|^\alpha \;
		    \|_{L_{\tilde{q}}(\rz\times X)},
	    \end{array}
	\end{equation}
	where $\tilde{\gamma}$ is defined by
	$\alpha \tilde{\gamma} - (\alpha-1) \sfrac{n+1}{2} = 
	\gamma + \beta$, and the first inequality holds in view of Lemma 
	\ref{badembed}. 

   In case $\delta<\frac{n+1}{2}$, we can decrease $\beta$ and assume 
   $\alpha \le \sfrac{n+1-2\gamma-2\beta}{n+1-2\delta}$. This implies 
   $\tilde{\gamma} \le \delta$. For $\delta\ge\frac{n+1}{2}$, note that 
   $\tilde{\gamma}=\frac{n+1}{2}(1-\frac{1}{\alpha})+\frac{\gamma+\beta}{\alpha}$. 
   Possibly decreasing $\beta$ , we obtain 
   $\tilde{\gamma}<\frac{n+1}{2}(1-\frac{1}{\alpha})+\frac{\delta}{\alpha}$, hence again 
   $\tilde{\gamma}- \delta\le 0$. But then
	\[
	\| S_{\tilde{\gamma}} u \|_{L_{\alpha \tilde{q}}(\rz \times X)}
	=
	\| \; |u|^\alpha \; 
	\|^\frac{1}{\alpha}_{\calH^{0,\gamma+\beta}_{\tilde{q}}(X^\wedge)}
	\le
	C \| u \|_{\calH^{s,\tilde{\gamma}}_{q}(X^\wedge)}
	\le 
	C \| u \|_{\calH^{s,\delta}_{q}(X^\wedge)},
	\]
	where the first inequality holds by Lemma \ref{contF} (with 
	$\tilde{\varrho} = \gamma + \beta$ and $\varrho = \tilde{\gamma}$)
	and the second one is true due to the boundedness of $\supp u$.
	Hence, $S_{\tilde{\gamma}}u$ runs through a bounded set of
	$L_{\alpha \tilde{q}}(\rz \times X)$ if $u$ runs through a bounded 
	set of $\calH^{s,\delta}_{q}(X^\wedge)_{0}$. We employ once more the
observation that the map 
	\[
	v \mapsto |v|^\alpha \; : \; L_{\alpha \tilde{q}}(\rz \times X)
	\rightarrow L_{\tilde{q}}(\rz \times X)
	\]
	is Lipschitz continuous on bounded sets and obtain from
\eqref{eq:lip1}
	\[
	\| \; |u|^\alpha - |v|^\alpha \; 
	\|_{\calH^{0,\gamma}_{q}(X^\wedge)}
	\le
	C \| S_{\tilde{\gamma}}u - S_{\tilde{\gamma}}v \|_{L_{\alpha 
	\tilde{q}}(\rz \times X)}
	\]
	for $u$ and $v$ in a bounded set of 
	$\calH^{s,\delta}_{q}(X^\wedge)_{0}$. Using the embeddings
	$H^s_{q}(\rz\times X) \hookrightarrow L_{\alpha \tilde{q}}(\rz\times 
	X)$ and $\calH^{s,\delta}_{q}(X^\wedge)_{0} \hookrightarrow
	\calH^{s,\tilde{\gamma}}_{q}(X^\wedge)_{0}$, we arrive at
	\[
	\| \; |u|^\alpha - |v|^\alpha \; 
	\|_{\calH^{0,\gamma}_{q}(X^\wedge)}
	\le
	C \| S_{\tilde{\gamma}}u - S_{\tilde{\gamma}}v 
	  \|_{H^{s}_{q}(\rz \times X)}
	=   C \| u - v \|_{\calH^{s,\tilde{\gamma}}_{q}(X^\wedge)}
	\le C \| u - v \|_{\calH^{s,\delta}_{q}(X^\wedge)},
	\]
	showing the desired Lipschitz continuity of \eqref{lip}.
    \end{proof}
    
    We are now ready to prove Theorem 
    \ref{nonlinear}. We know from Corollary
\ref{interp2} that 
    $E_{1-\frac{1}{q},q}\hookrightarrow\calH^{s,\delta}_q(\bz)$ for all
    $0<s<2/q^\prime$ and $\delta=\gamma_p+{2}/{q^\prime}-\varepsilon$
    with arbitrarily small $\varepsilon>0$. 
    Theorem \ref{lipF2} tells us when
$f:\calH^{s,\delta}_q(\bz)\to\calH^{0,\gamma_p}_q(\bz)$ 
is 
    Lipschitz continuous:
    In case ${2p}/{q^\prime}\le n+1$, we  have $\delta<({n+1})/{2}$;
    hence \eqref{condalpha1} -- or its simplified version, in case $sq \ge 
    n+1$ -- gives the admitted range of $\alpha$.
Similarly, ${2p}/{q^\prime}> n+1$ allows us to choose
$\delta>({n+1})/2$; this leads to 
\eqref{condalpha2} -- with the corresponding simplification for $sq\ge n+1$. 
Inserting $\gamma_p$ and $\delta$ in the expressions, letting $\eps \to 0$,
and optimizing 
$0<s<2/q'$, we obtain the formula for $\alpha^*$ in \ref{nonlinear}a) 
for $2p/q'<n+1$ and that in \ref{nonlinear}b) for $2p/q'\ge n+1$. 
Note that $\frac{2}{q^\prime}<\frac{n+1}{q}$ is equivalent to $q<\frac{n+3}{2}$.

\subsection{Conclusion} In order to illustrate the results, let us state one of
the possible applications of Proposition \ref{bDelta}, Theorem \ref{aDelta},
and Theorem \ref{nonlinear}.  
Others can be made up easily using Remark
\ref{integer},
Example \ref{experiodic}, and Remark \ref{u^alpha}. 
As before, 
$E_0=\calH^{0,\gamma_{p}}_{q}(\bz)$ and $
E_1=\calH^{2,\gamma_{p}+2}_{q}(\bz)$ with $1 < p,q < \infty$.

\begin{theorem} Let $n\ge 4$. 
Given $c>0$, $\alpha \ge 1$,  and
$T_0>0$,  
there is a suitable choice of $p$ and $q$ in
$]1,\infty[$ and $T_1>0$  such that the equation 
$$
\dot u - a(t^c u)\Delta u = |u|^\alpha +g \ \ \text{on }]0,T_0[, \ \ \ 
u(0)= u_0,
$$ 
has a unique solution $u\in  W^1_q([0,T_{1}];E_{0}) \cap L_q([0,T_{1}];E_{1})
  \cap {\mathcal C}([0,T_1],E_{1-\frac1q,q})$ on $]0,T_1[$ 
 for every $g\in L_q(]0,T_0[,E_0)$, every
smooth, strictly positive function $a$,
and every  $u_0\in \calC^\infty_{\rm comp}({\rm int}\,\bz)$.

Applying Remark {\rm \ref{integer}} we see that, 
for suitable $p,q\in{]1,\infty[}$, also the Ginzburg-Landau type equation
$$\dot u -\Delta u = u-u^3\ \ \text{on }]0,T_0[, \ \ \ u(0) = u_0,$$
has a unique solution $u$ on $]0,T_1[$ in the same space as above for
an arbitrary initial value $u_0$ in $E_{1-\frac1q,q}$.
\end{theorem}


\setlength{\parskip}{0pt}

\begin{small}
\bibliographystyle{amsalpha}

\end{small}

\end{document}